\documentclass[12pt,draftcls,onecolumn]{IEEEtran}

\usepackage{amssymb,amsbsy,amsmath,amsfonts,amscd,mathrsfs}
\usepackage{graphicx,verbatim}
\usepackage{color}

\newcommand{\llabel}[1]{{\label{#1}}}

\newcommand{\ffoot}[1]{}

\newcommand{\bi}{\begin{itemize}}
\newcommand{\ei}{\end{itemize}}
\newcommand{\bd}{\begin{description}}
\newcommand{\ed}{\end{description}}
\newcommand{\be}{\begin{enumerate}}
\newcommand{\ee}{\end{enumerate}}

\newcommand{\beq}{\begin{equation}}
\newcommand{\eeq}{\end{equation}}
\newcommand{\bqn}{\begin{eqnarray}}
\newcommand{\eqn}{\end{eqnarray}}
\newcommand{\eqnn}{\nonumber\end{eqnarray}}
\newcommand{\eqnl}[1]{\llabel{#1}\end{eqnarray}}

\newcommand{\nn}{\nonumber}

\newcommand{\ba}[1]{\begin{array}{#1}}
\newcommand{\ea}{\end{array}}

\newcommand{\R}{\mathbb{R}}
\newcommand{\C}{\mathbb{C}}
\newcommand{\N}{\mathbb{N}}

\newcommand{\fine}{\end{document}}

\def \trait (#1) (#2) (#3){\vrule width #1pt height #2pt depth #3pt}
\def \qed{\hfill
        \trait (0.1) (6) (0)
        \trait (6) (0.1) (0)
        \kern-6pt   
        \trait (6) (6) (-5.9)
        \trait (0.1) (6) (0)
\medskip}

\newtheorem{Theorem}{\bf Theorem}[section]
\newtheorem{ml}[Theorem]{\bf Lemma}
\newtheorem{example}[Theorem]{\bf Example}

\newtheorem{remark}[Theorem]{Remark}
\newtheorem{mcc}[Theorem]{\bf Corollary}
\newtheorem{Definition}[Theorem]{\bf Definition}
\newtheorem{mpr}[Theorem]{\bf Proposition}
\newtheorem{mproperty}[Theorem]{\bf Property}

\newcommand{\bt}{\begin{Theorem}}
\newcommand{\et}{\end{Theorem}}
\newcommand{\bl}{\begin{ml}}
\newcommand{\el}{\end{ml}}
\newcommand{\bp}{\begin{mpr}}
\newcommand{\ep}{\end{mpr}}
\newcommand{\bcc}{\begin{mcc}}
\newcommand{\ecc}{\end{mcc}}
\newcommand{\bex}{\begin{example}}
\newcommand{\eex}{\end{example}}

\newcommand{\bproperty}{\begin{mproperty}}
\newcommand{\eproperty}{\end{mproperty}}

\newcommand{\ec}{\end{mcc}}
\newcommand{\bdeff}{\begin{Definition}}
\newcommand{\edeff}{\end{Definition}}

\newcommand{\brem}{\begin{remark}}
\newcommand{\erem}{\end{remark}}

\newcommand{\proof}{\noindent {\bf Proof. }}

\newcommand{\lam}{\lambda}

\newcommand{\al}{\alpha}
\newcommand{\eps}{\varepsilon}

\newcommand{\con}{{\cal C}}

\renewcommand{\H}{{\cal H}}

\newcommand{\eproof}{\hfill $\Box$ \medskip}

\newcommand{\la}{\left\langle}
\newcommand{\ra}{\right\rangle}

\newcommand{\h}{h}

\newcommand{\curvepao}{Non-mixing curves}
\newcommand{\mariospectrum}{separated discrete spectrum}
\newcommand{\mariomatrix}{conicity matrix}
\newcommand{\campopao}{non-mixing field}
\newcommand{\Xpao}{\mathcal{X}_P}

\newcommand{\Cspet}{\con^{\infty}}

\newcommand{\uu}{{\bf u}}
\newcommand{\xx}{{\bf x}}
\newcommand{\Pdue}{\mathbb{P}}
\newcommand{\Hdue}{\mathbb{H}}
\newcommand{\Heff}{H^{\varepsilon}_{\mathrm{eff}}}
\newcommand{\Ueff}{U^{\varepsilon}_{\mathrm{eff}}}
\newcommand{\tphi}{\widetilde{\phi}}

\newcommand{\tg}{\widehat{\gamma}}

\newcommand{\tro}{\widehat{\rho}}
\newcommand{\tet}{\widehat{\theta}}

\newcommand{\bth}{\boldsymbol{\Xi}}

\newcommand{\spazio}{\con^0(\Omega,\R)}

\title{
Adiabatic control of the Schr\"odinger equation via conical intersections of the eigenvalues
}

\author{U.~Boscain \thanks{CMAP \'Ecole Polytechniques CNRS, Palaiseau, France \texttt{boscain@cmap.polytechnique.fr
}} \and F. Chittaro \thanks{LSS-Sup\'elec, 3 rue Joliot-Curie, 91192 Gif-sur-Yvette, France \texttt{francesca.chittaro@lss.supelec.fr}} \and P.~Mason \thanks{CNRS-LSS-Sup\'elec, 3 rue Joliot-Curie, 91192 Gif-sur-Yvette, France \texttt{paolo.mason@lss.supelec.fr}} \and M.~Sigalotti \thanks{
CORIDA, INRIA Nancy -- Grand Est, France and 
 Institut
\'Elie Cartan, UMR 7502 Nancy-Universit\'e/CNRS,
BP 239, Vand\oe uvre-l\`es-Nancy 54506,  France 
\texttt{mario.sigalotti@inria.fr}} 
\and \thanks{This research has been supported  by the European Research Council, ERC
StG 2009 ``GeCoMethods", contract number 239748, by the ANR ``GCM", program ``Blanc--CSD"
project number NT09-504490, and by the DIGITEO project ``CONGEO".
}
}

\begin{document}
\maketitle

\begin{abstract}
In this paper we present a constructive method to control the bilinear Schr\"odinger equation via two controls.
The method is based on adiabatic techniques and works if the spectrum
of the Hamiltonian admits eigenvalue intersections, and if the latter
are conical (as it happens generically).
We provide sharp estimates of the relation between the error and the controllability time. 
\end{abstract}

\section{Introduction}
In this paper we are concerned with the problem of controlling the Schr\"odinger equation
\bqn
i\frac{d\psi}{dt}=\left(H_0+\sum_{k=1}^m u_k(t)H_k\right)\psi(t).
\label{eqMA}
\eqn
Here $\psi$ belongs to the Hilbert sphere  $\mathbf{S}$ of a complex separable Hilbert space ${\cal H}$ and $H_0,\ldots, H_m$  
are self-adjoint operators on ${\cal H}$.
The controls $u_1,\ldots,u_m$ are scalar-valued and represent the action of  external fields. $H_0$ describes the ``internal'' dynamics of the system, while $H_1,\ldots, H_m$  the interrelation between the system and the controls.

The reference model is the one in which  $H_0=-\Delta+V_0(x)$, $H_i=V_i(x)$, where $x$ belongs to a domain $D\subset\R^n$  and $V_0,\dots,V_m$ are real functions (identified with the corresponding multiplicative operators). 
However, equation \eqref{eqMA} can be used to describe  more general controlled dynamics. For instance, a quantum particle on a Riemannian manifold subject to  external fields (in this case $\Delta$ is the  Laplace--Beltrami operator) 
or a two-level ion trapped in a harmonic potential (the so-called Eberly and Law model~\cite{Boscain_Adami,rangan}).
In the last case, as in many other relevant physical situations, the operator $H_0$ cannot be written as the sum of a Laplacian plus a potential.

The controllability problem consists in establishing whether, for every pair of states $\psi_0$ and $\psi_1$, there exist  controls $u_k(\cdot)$ and a time $T$ such that the solution of~\eqref{eqMA} with initial condition $\psi(0)=\psi_0$ satisfies $\psi(T)=\psi_1$. The answer to this question is negative when ${\cal H}$ is infinite-dimensional. Indeed,  Ball, Marsden and Slemrod proved in~\cite{bms} a result which implies (see~\cite{turinici}) that equation \eqref{eqMA} is not controllable in (the Hilbert sphere of) ${\cal H}$. 
Hence one has to look for weaker controllability properties as, for instance, approximate controllability  or controllability  between 
the eigenstates of $H_0$ (which are the most relevant physical states).
However, in certain cases one can describe quite precisely the set of states which can be connected by admissible paths (see \cite{beauchard-coron,camillo,fratelli-nersesyan}).

In~\cite{schNOI1} an approximate controllability result for~\eqref{eqMA} was proved via finite-dimensional geometric control techniques applied to the Galerkin approximations.  
The main hypothesis is that the spectrum of $H_0$ is discrete and without rational resonances, which means that the gaps between the eigenvalues of $H_0$ should be $\mathbb{Q}$-linearly independent. Another crucial hypothesis is that the operator $H_1$ couples all eigenvectors of $H_0$. This result has been improved in \cite{schNOI2} where the hypothesis of 
 $\mathbb{Q}$-linear independence was weakened.
Similar results 
have been obtained, with different techniques, in \cite{nersesyan} (see also \cite{mirrahimi-continuous,fratelli-nersesyan}).

The practical application of the results discussed above entails three main difficulties:

$\bullet$ In most cases the techniques used to get controllability results do not permit to obtain (even numerically) the controls necessary to steer the system between two given states. 

$\bullet$ Even in the cases in which  one can get the controls as a byproduct of the controllability result, they happen to be  highly oscillating and hence they  can be difficult to implement, depending on the experimental conditions. 
Roughly speaking, since one  should move in an infinite dimensional space with only one control, one should generate many iterated Lie brackets. This is particularly evident in the papers \cite{schNOI2,schNOI1}, where the use of Galerkin approximations permits to 
highlight the Lie algebra structure.

$\bullet$ Explicit expressions of time estimates, for the norm of controls and for their total variations are extremely difficult to obtain. (For some lower bounds of controllability time and estimates of the $L^1$ norm of the controls see \cite{schNOI2}.)

In most of the results in the literature  only the case $m=1$  is considered. 
In this paper we study the case $m=2$ and we  get both controllability results and explicit expressions of the external fields realizing the transition. The system under consideration is then
\[i\frac{d}{dt}\psi(t)=H(u_1(t),u_2(t))\psi(t),\]  
with $H(u_1,u_2)=H_0+u_1H_1+u_2H_2$.
The idea is to use  two slowly varying controls and climb the energy levels through conical intersections, if they are present. Conical eigenvalue intersections have been used to get population transfers in the finite dimensional case in \cite{u-STIRA,jauslin2,rouchon_a,ticozzi,jauslin}. Some preliminary ideas given in the present paper can be found in \cite{Boscain_Adami}, where a specific example (which is  a version of the Eberly and Law model) is analyzed, and in \cite{cdc-atalanta}. The main ingredients of our approach are the following:

$\bullet$ The adiabatic theorem that, in its rougher form, states the following:  
let $\lam(u_1,u_2)$ be an eigenvalue of $H(u_1,u_2)$ depending continuously on $(u_1,u_2)$ and assume that, for every $u_1,u_2\in K$ ($K$ compact subset of $\R^2$), $\lam(u_1,u_2)$ is simple.  Let $\phi(u_1,u_2)$ be the corresponding eigenvector (defined up to a phase). Consider a path $(u_1,u_2):[0,1]\to K$  and its reparametrization  $(u_1^\eps(t),u_2^\eps(t))=(u_1(\eps t),u_2(\eps t))$, defined on $[0,1/\eps]$. Then the solution $\psi_\eps(t)$ of the equation  
$i\frac{d\psi_\eps}{d  t }=(H_0+u_1^\eps (t)H_1+u_2^\eps (t)H_2)\psi_\eps(t)$ with initial condition $\psi_\eps(0)=\phi(u_1(0),u_2(0))$ satisfies
  \bqn
\left\|\psi_{\eps}\left({1}/{\eps}\right)-e^{i \vartheta}\phi\left(u_1^\eps\left( 1 /\eps\right),u_2^\eps\left( 1/ \eps\right)\right)\right\|\leq C\eps
 \label{eq-CgapSQUARE}
 \eqn
 for some $\vartheta=\vartheta(\eps)\in\R$.  This means that, if the controls are slow enough, then, up to phases, the state of the system follows the evolution of the eigenstates of the time-dependent Hamiltonian. The constant $C$ depends on the gap between the eigenvalue $\lam$ and the other eigenvalues. 

\begin{figure}
\centering
\input{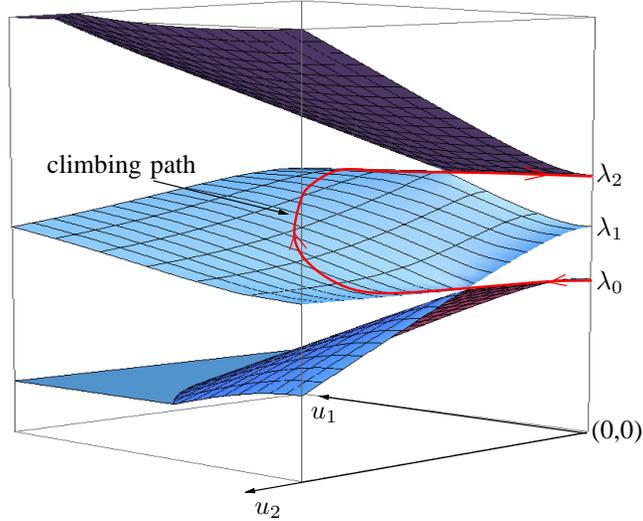}
\caption{A slow path climbing the spectrum of $H(u_1,u_2)$, plotted in function of $(u_1,u_2)$.}
\label{s-POST-climbing}
\end{figure}

$\bullet$ The crossing of conical intersections. 
Generalizations of the adiabatic theory guarantee that, if the path $(u_1(\cdot),u_2(\cdot))$
 passes (once) through a conical intersection between the eigenvalues $\lambda_0\leq\lambda_1$,  
then   
\bqn
 \|\psi_{\eps}(1/\eps)-e^{i \vartheta}\phi_1(u_1^\eps(1/\eps),u_2^\eps(1/\eps))\|\leq C\sqrt{\eps}
 \label{eq-Cgap}
 \eqn
where $\psi_\eps(t)$ is the solution of the equation  
$i\frac{d\psi_\eps}{d  t }=(H_0+u_1^\eps (t)H_1+u_2^\eps (t)H_2)\psi_\eps(t)$ with initial condition $\psi_\eps(0)=\phi_0(u_1^\eps(0),u_2^\eps(0))$ and $\phi_0$, $\phi_1$ are the eigenvectors corresponding respectively to the eigenvalues $\lam_0$, $\lam_1$ (see \cite{teufel}).
Figure~\ref{s-POST-climbing} illustrates a closed slow path in the space of controls producing a transition from the eigenvector corresponding to the eigenvalue $\lam_0$ to the eigenvector corresponding to $\lam_2$ by crossing two conical singularities. Notice that the path could not be closed if only one control was present. Indeed $u(t)$ would pass back and forth through the same singularity and the trajectory would come back to the original state. 
One of our main results is that choosing special curves that pass through the conical singularity the estimate in \eqref{eq-Cgap} can be improved by replacing $\sqrt{\eps}$ by $\eps$. 
Hence if some energy levels $\lam_0,\ldots, \lam_k$ of the spectrum of $H$ are connected by conical singularities, then one can steer, in time $1/\eps$, an eigenstate corresponding to $\lam_0$ to an eigenstate corresponding to $\lam_k$  with an error of order ${\eps}$. 

$\bullet$ The behavior of the eigenstates in a neighborhood of conical singularities. If the path $(u_1(\cdot),u_2(\cdot))$ is a piecewise smooth curve with a \emph{vertex} (i.e.  discontinuity at the $\con^1$ level) at the conical singularity,
 the state of the system evolves with continuity,
while the eigenstates corresponding to the degenerate eigenvalue are subject to an instantaneous rotation. The angle made by the path at the vertex can be used to control the splitting of probabilities between the two energy levels (see Figure~\ref{s-POST-climbing}). This splitting phenomenon has already been described and exploited for controllability purposes on a two-dimensional system  in~\cite{cdc-atalanta}.

\begin{figure}
\begin{center}
\input{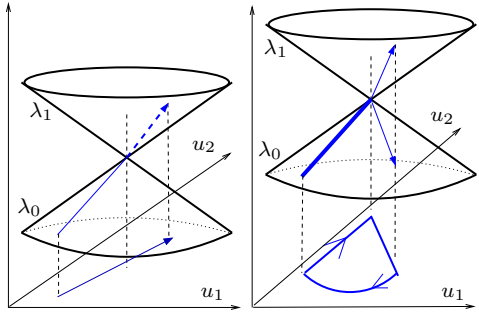}
\caption{Passages through a conical intersection.}
\label{s-angolino}
\end{center}
\end{figure}

The ideas introduced above lead to the following result:
  if the  energy levels $\lam_0,\ldots, \lam_m$ are connected 
by conical singularities, then the system is approximately {\sl spread controllable}, i.e.,  for every given $\eps>0$ and $p_0,\ldots,p_m\ge0$ such that  $\sum_{i=0}^m p_i^2=1$, there exists a control ${\bf u}$ defined on $[0,1/\eps]$,  $k+1$ phases $\vartheta_0,\ldots,\vartheta_k \in \mathbb{R}$, and a trajectory corresponding  to ${\bf u}$ satisfying  $\psi(0)=\phi_0$ and   $\|\psi(1/\eps)-\sum_{j=1}^k p_j e^{i\vartheta_j} \phi_j(\uu^1)\|\leq \varepsilon.$
Moreover the control can be taken of the form
${\bf u}(t)=\gamma(\eps t)$,  where $\gamma:[0,1]\to \R^2$ is characterized explicitly. 
 Hence the method provides precise time estimates in relation with the required precision.
The method cannot be easily reversed, in order to explicitly characterize  paths 
steering a state which is spread on several eigenstates to a single one. 
The difficulty lies on the loss of information about the relative phases 
during adiabatic evolution. 

We finally remark that systems for which the method can be applied are rather frequent. Indeed intersections of eigenvalues are generically conical for Hamiltonians of the form $-\Delta+V_0+u_1 V_1+u_2 V_2$, as explained in Section \ref{definizionieris}.

The structure of the paper is the following.
In Section~\ref{definizionieris}, we introduce the framework and we state the main result. In Section~\ref{literature} we recall the time adiabatic theorem and some results on the regularity of eigenvalues and eigenstates of parameter-dependent Hamiltonians.
In Section~\ref{sconical} we deepen our analysis of conical intersection; in particular, we state and prove a sufficient condition for an intersection to be conical.
Section~\ref{una-a-caso} is devoted to the construction of some special curves along which we can obtain our controllability result, while the proof of the main theorem is the subject of Section~\ref{prova}.
Finally, in Section~\ref{mild}, we show that the same controllability result holds also for more general curves than those presented in Section~\ref{una-a-caso}.


\section{Definitions and main results}
\label{definizionieris}

We consider the Hamiltonian 
\[
H(\uu) = H_0 + u_1 H_1 +u_2 H_2,
\]
for $\uu = (u_1, u_2) \in \R^2$. From now on we assume that $H(\cdot)$ satisfies the following assumption:

\smallskip
\noindent
{\bf (H0)} $H_0$ is a self-adjoint operator on a separable Hilbert space ${\cal H}$, and $H_1$ and $H_2$ are bounded self-adjoint operators on ${\cal H}$.

\smallskip

\noindent When necessary, we also make the following assumption on the Hamiltonian $H(\cdot)$:

\smallskip

\noindent
{\bf (H1)} There exists an orthonormal basis $\{\chi_j\}_j$ of the Hilbert space $\mathcal{H}$ such that the matrix elements $\langle \chi_j, H_0 \chi_k \rangle$, $\langle \chi_j, H_1 \chi_k \rangle$ and $\langle \chi_j, H_2 \chi_k \rangle$ are real for any $j,k$.

\smallskip

\brem 
Hypothesis {\bf (H1)} ensures that, with each $\uu$ and each eigenvalue of $H(\uu)$ (counted according to their multiplicity), it is possible to associate an eigenstate whose components with respect to the basis $\{\chi_j\}_j$ are all real.\label{reale}
\erem

A typical case for which {\bf (H0)} and {\bf (H1)} are satisfied is when $H_0=-\Delta  + V$, where $\Delta$ is the Laplacian on a bounded domain $\Omega\subset \R^d$ with Dirichlet boundary conditions,  $V\in L^\infty(\Omega,\R)$, ${\cal H}=L^2(\Omega,\C)$, 
and $H_1, H_2$ are two bounded multiplication operators by real valued functions. In this case the spectrum of $H_0$ is  discrete. 

The dynamics are described by the time-dependent Schr\"odinger equation
\bqn
i\frac{d\psi}{dt}=H(\uu(t))\psi(t).
\label{eqGRA}
\eqn
Such equation has classical solutions under hypothesis {\bf (H0)}, $\uu(\cdot)$ piecewise $\con^1$ and with an initial condition in the domain of $H_0$ (see \cite{segal63} and also \cite{bms}).

We are interested in controlling  \eqref{eqGRA} inside some portion of the discrete spectrum of $H(\uu)$. 
Since we use adiabatic techniques, the structure of the spectrum shall satisfy some particular features: roughly, the portion of discrete spectrum we consider must be well separated from its complement in the spectrum of the Hamiltonian, and this property must hold uniformly for $\uu$ belonging to some domain in $\R^2$.

All these properties are formalized by the following notion:

\bdeff 
Let $\omega$ be a domain in $\R^2$. A map $\Sigma$ defined on $\omega$ that associates with each $\uu\in\omega$ a subset  $\Sigma(\uu)$ of the discrete spectrum of $H(\uu)$ is said to be a \emph{\mariospectrum}  on $\omega$ if there exist two continuous functions $f_1,f_2:\omega\to\R$ such that:
\bi
\item $f_1(\uu) < f_2(\uu)$ and $\Sigma(\uu)\subset [f_1(\uu) , f_2(\uu)]\qquad\forall \uu\in\omega$.
\item there exists $\Gamma>0$ such that
\[\inf_{\uu\in\omega} \inf_{\lambda \in \mathrm{Spec}(H(\uu))\setminus \Sigma(\uu)} \mathrm{dist}(\lambda,[f_1(\uu), f_2(\uu)])  )>\Gamma.
\]
\ei
\edeff 


\noindent{\bf Notation}
From now on we label  the eigenvalues  belonging to $\Sigma(\uu)$ in such a way that we can write $\Sigma(\uu)=\{\lambda_0(\uu) ,\ldots,\lambda_k(\uu)\}$, where  
$\lambda_0(\uu)\leq \cdots \leq\lambda_k(\uu)$ are counted according to their multiplicity (note that the separation of $\Sigma$ from the rest of the spectrum guarantees that $k$ is constant).
Moreover we denote by  $\phi_0(\uu) , \ldots , \phi_k(\uu)$ an orthonormal family of eigenstates corresponding to $\lambda_0(\uu), \ldots , \lambda_k(\uu)$.  Notice that in this notation $\lambda_0$ needs not being the ground state of the system.\\[2mm]

\bdeff
Let $\Sigma$ be a \mariospectrum\ on $\omega$.
 We say that \eqref{eqGRA} is approximately \emph{spread-controllable} on $\Sigma$ if for every 
 $\uu^0,\uu^1\in \omega$ such that 
 $\Sigma(\uu^0)$ and $\Sigma(\uu^1)$ are non-degenerate, 
for every 
 $\bar{\phi}\in \{\phi_0(\uu^0) ,\ldots,\phi_k(\uu^0) \}$, $p \in [0,1]^{k+1}$ such that $\sum_{l=0}^k p_l^2=1$,  and every $\eps>0$ there exist $T>0$, $\vartheta_0,\ldots,\vartheta_k \in \mathbb{R}$ and a piecewise $\con^1$ control $\uu(\cdot):[0,T]\to \R^2$
such that 
\begin{equation} \label{app-ctr}
\|\psi(T)-\sum_{j=0}^k p_j e^{i\vartheta_j} \phi_j(\uu^1)\|\leq \varepsilon,
\end{equation}
where $\psi(\cdot)$ is the solution of \eqref{eqGRA} with $\psi(0)=\bar{\phi}$.  \label{delfozorzi}
\edeff


Our techniques rely on the existence of conical intersections between the eigenvalues. 
Conical intersections constitute a well-known notion in molecular physics. 
They have an important role in the Born--Oppenheimer approximations (see for instance \cite{bofo,lasser,teufel}, where they appear for finite dimensional operators). In the finite dimensional case they have been classified by Hagedorn \cite{hagedorn}. 

A unified 
characterization 
of conical intersections
seems to be missing. The following definition meets
 all the features commonly 
 attributed to 
 them. 

\bdeff
Let $H(\cdot)$ satisfy hypothesis {\bf (H0)}. We say that $\bar\uu\in\R^2$ is a \emph{conical intersection} between the eigenvalues $\lam_j$ and $\lam_{j+1}$ if $\lam_{j}(\bar \uu) = \lam_{j+1}(\bar \uu)$ has multiplicity two and there exists a constant $c>0$ such that for any unit vector $\mathbf{v}\in \R^2$ and $t>0$ small enough we have that 
\begin{equation} \label{formcono}
\lam_{j+1}(\bar \uu+t\mathbf{v})-\lam_{j}(\bar \uu+t\mathbf{v}) > ct\,.\end{equation}
\label{conical}
\edeff

It is worth noticing that conical intersections are not pathological
phenomena. On the contrary, they happen to be generic
in the following sense.
Consider the reference case where $\H=L^2(\Omega,\C)$,
$H_0=-\Delta+V_0:D(H_0)=H^2(\Omega,\C)\cap H^1_0(\Omega,\C)\to
L^2(\Omega,\C)$, $H_1=V_1$, $H_2=V_2$, with $\Omega$ a bounded domain
of $\R^d$ for some $d\in\N$ and $V_j\in \spazio$ for $j=0,1,2$.
Then, generically with respect to the pair $(V_1,V_2)$ in
$\spazio\times \spazio$ (that is, for all pairs, $(V_1,V_2)$ in a
countable intersection of open and dense  subsets of $\spazio\times \spazio$),
for each $\uu \in \mathbb{R}^2$ and $\lambda\in \R$ such that $\lambda$ is a multiple
eigenvalue of
$H_0+u_1 H_1+u_2 H_2$, the eigenvalue intersection $\uu$ is conical.

In order to check that this is true,
we can
apply the transversal density theorem (see \cite[Theorem
19.1]{Abraham-Robbin}) with $\mathscr{A}=\spazio\times \spazio$,
$X=\R^2$,
$Y= \spazio$,
$\rho((V_1,V_2),\uu)=V_0+u_1 V_1+u_2 V_2$,
and
$$W=\{V\in \spazio\mid -\Delta+V:H^2(\Omega,\C)\cap
H^1_0(\Omega,\C)\to L^2(\Omega,\C)\mbox{ has multiple
eigenvalues}\}.$$
The covering of $W$ by manifolds of codimension two is obtained in
\cite{teytel}, based on the properties proved in
\cite{besson-commentari-89} (see also \cite{lambert-lanzaMPAG2006}).
We obtain that, generically with respect to $(V_1,V_2)$,
the intersection of $\rho((V_1,V_2),\R^2)$ with $W$ is transverse.
Equivalently said, generically with respect to $(V_1,V_2)$,
for every $\uu \in\mathbb{R}^2$ and $\lambda\in \R$ such that $\lambda$ is a multiple
eigenvalue of
$-\Delta+V_0+u_1 V_1+u_2 V_2$, for every
$(v_1,v_2)\in\R^2\setminus\{0\}$,  the line $\{(u_1+tv_1)
V_1+(u_2+tv_2) V_2\mid t\in\R\}$ is not tangent to $W$, i.e., the
eigenvalue intersection $\uu$ is conical.

Moreover, each
conical intersection $(u_1,u_2)$ is structurally stable, in the sense
that  small perturbations of $V_0$, $V_1$ and $V_2$ give rise, in a neighborhood of $\uu$, to conical
intersections for the perturbed $H$.
Structural stability properties can be proved without resorting to
abstract transversality theory, as
will be shown in Section~\ref{una-a-caso}, Theorem~\ref{strutto}.

Our main result is the following: it states that spread controllability holds for a class of systems having pairwise conical intersections, providing in addition an estimate of the controllability time. As a byproduct of the proof, we will also get an explicit characterization of the motion planning strategy (the path $\gamma(\cdot)$ below).
\bt
\label{asc-paolo} 
Let $H(\uu)=H_0+u_1 H_1+u_2 H_2$ satisfy hypotheses {\bf (H0)-(H1)}.
Let $\Sigma :\uu \mapsto \{\lam_0(\uu),\ldots,\lam_k(\uu) \}$  be a \mariospectrum\ on $\omega\subset \R^2$ and assume that there exist conical intersections $\uu_j\in\omega,\ j=0,\ldots,k-1$, between the eigenvalues $\lam_j,\lam_{j+1}$, with $\lambda_l(\uu_j)$ simple if $l\neq j,j+1$.
Then, 
for every $\uu^0$ and $\uu^1$ such that $\Sigma(\uu^0)$ and $\Sigma(\uu^1)$ are non-degenerate, for every
$\bar{\phi}\in \{\phi_0(\uu^0) ,\ldots,\phi_k(\uu^0) \}$, and $p \in [0,1]^{k+1}$ such that $\sum_{l=0}^k p_l^2=1$, 
there exist $C>0$ and a continuous control $\gamma(\cdot):[0,1]\to \R^2$ with $\gamma(0)=\uu^0$ and $\gamma(1)=\uu^1$, such that for every $\eps>0$
\[
\|\psi(1/\eps)-\sum_{j=0}^k p_j e^{i\vartheta_j} \phi_j(\uu^1)\|\leq C \varepsilon,
\]
where $\psi(\cdot)$ is the solution of (\ref{eqGRA}) with $\psi(0)=\bar{\phi}$, $\uu(t)=\gamma(\eps t)$, and $\vartheta_0,\ldots,\vartheta_k \in \mathbb{R}$ are some phases depending on $\eps$ and $\gamma$.  
In particular, 
(\ref{eqGRA}) is approximately spread controllable on $\Sigma$. 
\et

\section{Survey of basic results}\label{literature}

\subsection{The adiabatic theorem} \label{subs ad}

One of the main tools used in this paper is the adiabatic theorem (\cite{bofo,kato,nenciu,panati}); here we recall its formulation, adapting it to our framework. For a general overview see the monograph \cite{teufel}. 
We remark that we refer here exclusively to the time-adiabatic theorem.

The adiabatic theorem deals with quantum systems governed by Hamiltonians that explicitly depend on time, but whose dependence is slow. While in quantum systems driven by time-independent Hamiltonians the evolution preserves the 
 occupation probabilities of the energy levels, this is in general not true for time-dependent Hamiltonians. The adiabatic theorem states that if the time-dependence is slow, then the occupation probability of the energy levels, which also evolve in time, is approximately conserved by the evolution.

More precisely, consider $\h(t)= H_0 + u_1(t) H_1 + u_2(t) H_2$, $t \in I=[t_0,t_f]$, 
satisfying {\bf (H0)}, and assume that the map $t \mapsto (u_1(t),u_2(t))$ belongs to $\con^2(I)$.
Assume moreover that there exists $\omega \subset \mathbb{R}^2$ such that $(u_1(t),u_2(t)) \in \omega$ for all $t \in I$ and $\Sigma$ is a \mariospectrum\ on $\omega$.

We introduce a small parameter $\eps>0$ that controls the time scale, and consider the slow Hamiltonian $\h(\eps t),$ $t \in [t_0/\eps,t_f/\eps]$.
The time evolution (from $t_0/\eps$ to $t$) $\widetilde{U}^{\eps}(t,t_0/\eps)$ generated by $\h(\eps \cdot)$ satisfies the equation
$i \frac{d}{dt}\widetilde{U}^{\eps}(t,t_0/\eps) = \h(\eps t) \widetilde{U}^{\eps}(t,t_0/\eps)$.
Let $\tau=\eps t$ belong to $[t_0,t_f]$ and $\tau_0 =t_0$;  the time evolution $U^{\eps}(\tau,\tau_0):=\widetilde{U}^{\eps}(\tau/\eps,\tau_0/\eps)$ satisfies the equation
\begin{equation} 
i \eps\frac{d}{d\tau}U^{\eps}(\tau,\tau_0) = \h(\tau) U^{\eps}(\tau,\tau_0). \\
\end{equation}
Notice that $U^{\eps}(\tau,\tau_0)$ does not preserve the probability of occupations: in fact, if we denote by $P_*(\tau)$ the spectral projection of $\h(\tau)$ on $\Sigma(\uu(\tau))$, then 
$P_*(\tau) U^{\eps}(\tau,\tau_0)$
is in general  
different from $U^{\eps}(\tau,\tau_0)P_*(\tau_0)$. 

Let us consider the \emph{adiabatic Hamiltonian} associated with $\Sigma$:
\[
\h_a(\tau) = \h(\tau) -i \eps P_*(\tau) \dot{P}_*(\tau) -i \eps P^{\bot}_*(\tau) \dot{P}^{\bot}_*(\tau),
\]
where $P^{\bot}_*(\tau)=\mathrm{id}-P_*(\tau)$ and $\mathrm{id}$ denotes the identity on $\mathcal{H}$. Here and in the following the time-derivatives 
 shall be intended with respect to the reparametrized time $\tau$.
The adiabatic propagator associated with $\h_a(\tau)$,  denoted by $U^{\eps}_a(\tau,\tau_0)$, is the solution of the equation 
\begin{equation} \label{adia-ev}
i \eps\frac{d}{d\tau}U_a^{\eps}(\tau,\tau_0) = \h_a(\tau) U_a^{\eps}(\tau,\tau_0)
\end{equation}
with $U_a^{\eps}(\tau_0,\tau_0)=\mathrm{id}$.

Notice that 
\[
P_*(\tau) U_a^{\eps}(\tau,\tau_0) = U_a^{\eps}(\tau,\tau_0)P_*(\tau_0), 
\]
that is, the 
adiabatic evolution preserves the occupation probability of the band $\Sigma$. 

Now we can  adapt to our setting the strong version of the quantum adiabatic theorem, as stated in \cite{teufel}.

\bt \label{th: adiabatic}
Assume that
$
{H}(\uu)= H_0 + u_1 H_1 + u_2 H_2$ satisfies {\bf (H0)}, and that $\Sigma$ is a \mariospectrum\ on $\omega \subset \mathbb{R}^2$. Let $I=[t_0,t_f]$,
$\uu : I \rightarrow \omega$ be a $\con^2$ curve and set $\h(t)={H}(\uu(t))$.
Then $P_* \in \con^2(I,\mathcal{L}(\mathcal{H}))$ and there exists a constant $C>0$ such that for all $\tau,\tau_0 \in I$
\begin{equation} \label{eq: adiabatic}
\|U^{\eps}(\tau,\tau_0)-U_a^{\eps}(\tau,\tau_0)\| \leq C\eps \left(1+ |\tau-\tau_0|\right).
\end{equation}
\et

\brem \label{nenciu}
If there are more than two parts of the spectrum which are separated by a gap, then it is possible to generalize the adiabatic Hamiltonian in the following way (\cite{nenciu}):
\[
\h_a(\tau) = \h(\tau) -i \eps \sum_{\alpha} P_{\alpha}(\tau) \dot{P}_{\alpha}(\tau) 
\]
where each $P_{\alpha}(\tau)$ is the spectral projection associated with a separated portion  of the spectrum, partitioning it as $\alpha$ varies.
\erem

\brem
In general the adiabatic theorem is stated for a time dependent
Hamiltonian $\h(t)$ satisfying the following hypotheses: 
it is assumed that all the Hamiltonians $\h(t)$ have a common dense domain $\mathcal{D}$ and that the function $t \mapsto \h(t)$ is $\con^2(I)$ and bounded as a function from $I$ to $\mathcal{L}_{sa}(\mathcal{D},\mathcal{H})$, where $\mathcal{L}_{sa}(\mathcal{D},\mathcal{H})$ denotes the space of bounded self-adjoint linear operators
from $\mathcal{D}$ to $\mathcal{H}$ and $\mathcal{D}$ is endowed with the norm of the graph of $\h(\bar{t})$, for some $\bar{t}\in I$. 
These hypotheses are satisfied for an Hamiltonian of the form
$\h(t)=H_0+u_1(t)H_1+u_2(t) H_2$ under assumption {\bf (H0)}, provided that the curve $\uu(\cdot)$ is $\con^2$. 
\erem

In this paper we are particularly interested in the adiabatic evolution across conical intersections between eigenvalues. A result in this direction can be found in
\cite[Corollary 2.5]{teufel}. In the language of control theory it reads as follows.
\bp\label{p_teufel}
Let $\Sigma: \uu \mapsto \{\lam_0(\uu),\ldots,\lam_k(\uu)\}$ be a \mariospectrum\ on $\omega$. Let $\uu^0,\uu^1,\bar\uu_j \in\omega,\ j=0,\ldots,k-1$. 
Assume that $\lambda_l(\uu^0)$ and $\lambda_l(\uu^1)$ are simple for all $l=0,\ldots,k$, and that, for any $j=0,\ldots,k-1$, $\bar\uu_j$ is a conical intersection between $\lambda_j$ and $\lambda_{j+1}$, with $\lambda_l(\uu_j)$ simple if $l\neq j,j+1$.
Let  $\gamma(\cdot):[0,1]\rightarrow \omega$ be a $\con^2$  curve
with $\gamma(0)=\uu^0$ and $\gamma(1)=\uu^1$ and such that 
 the eigenstates corresponding to  the eigenvalues $\lambda_l$  can be chosen $\con^1$ along $\gamma$ for all $l=0,\ldots,k$.
Assume moreover that  
 there exist times $0<\bar t_0<\cdots<\bar t_{k-1}<1$ with $\gamma(\bar t_j)=\bar\uu_j$, $\dot{\gamma}(\bar t_j)\neq 0, \ j=0,\ldots,k-1$, and that for any $l=1,\ldots,k\ $
$\lambda_l(\gamma(t))$ is
simple for every $t\neq \bar t_j,\ j=1,\ldots,k-1$. 
 
Then there exists $C>0$ such that, for any $\varepsilon>0$
\begin{equation} \label{teufel k}
\big\| \psi(1/\eps)-e^{i\vartheta }\phi_k\big\|\leq C \sqrt{\varepsilon}, 
\end{equation}
where $\vartheta \in \mathbb{R}$ and $\psi(\cdot)$ is the solution of equation 
\eqref{eqGRA} with $\psi(0)=\phi_0(\uu^0)$ corresponding to the control $\uu:[0,1/\eps]\rightarrow \omega$ defined by $\uu(t)=\gamma(\eps t)$.
\ep

In this paper we are interested in finding control paths along which we have a knowledge of adiabatic evolution finer than in \eqref{teufel k}. This allows also richer control strategies than those described in Proposition~\ref{p_teufel}, as it is needed to prove spread controllability. For this purpose we write an effective Hamiltonian describing the dynamics inside a two-dimensional band, possibly with conical intersections.

 Let us then consider the band constituted by the eigenvalues $\lambda_j,\lambda_{j+1} \in \Sigma$; we can find an open domain $\omega'\subset \omega$ such that $\{\lambda_j,\lambda_{j+1}\}$ is a \mariospectrum\ on $\omega'$. 

As above, we consider a control function $\uu(\cdot) \in \con^2(I,\omega')$, for a given time interval $I$. 
We can then apply the adiabatic theorem to the \mariospectrum\ $\Sigma': \uu \mapsto \{\lambda_j(\uu),\lambda_{j+1}(\uu)\}$,  $\uu \in \omega'$: we call $\Pdue(\tau)$ the spectral projection on the band $\{\lambda_j(\uu(\tau)),\lambda_{j+1}(\uu(\tau))\}$ and $\Hdue(\tau)=\Pdue(\tau)\mathcal{H}$ its range, which is the direct sum of the eigenspaces of $\lambda_j(\uu(\tau))$ and
$\lambda_{j+1}(\uu(\tau))$. 
We consider the adiabatic Hamiltonian $h_a(\tau) = h(\tau) -i \eps \Pdue(\tau) \dot{\Pdue}(\tau) -i \eps \Pdue^{\bot}(\tau) \dot{\Pdue}^{\bot}(\tau)$ and its associated propagator $U_a^{\eps}(\tau,\tau_0)$.

We are interested in describing the dynamics inside $\Hdue(\tau)$. Since  $\Hdue(\tau)$ is two-dimensional for any $\tau$, it is possible to map it isomorphically on $\mathbb{C}^2$ and identify an \emph{effective Hamiltonian} whose 
evolution is a representation of $U_a^{\eps}(\tau,\tau_0)|_{\Hdue(\tau_0)}$ on $\mathbb{C}^2$.

Let us assume that there exists an eigenstate basis $\{\phi_{\alpha}(\tau),\phi_{\beta}(\tau)\}$ of $\Hdue(\tau)$ such that 
$\phi_{\alpha}(\cdot),\phi_{\beta}(\cdot)$ belong to $\con^1(I,\mathcal{H})$.
We construct the time-dependent unitary operator $\mathcal{U}(\tau) : \Hdue(\tau) \rightarrow \mathbb{C}^2$ by defining for any $\psi \in \Hdue(\tau)$
\begin{equation} \label{Urep}
\mathcal{U}(\tau) \psi=e_1 \langle \phi_{\alpha}(\tau), \psi \rangle + e_2 \langle \phi_{\beta}(\tau), \psi \rangle,
\end{equation}
\noindent
where $\{e_1,e_2\}$ is the canonical basis of $\mathbb{C}^2$.

We then define the \emph{effective propagator}
\begin{equation} \label{effprop}
\Ueff(\tau,\tau_0) = \mathcal{U}(\tau)U_a^{\eps}(\tau,\tau_0)\mathcal{U}^*(\tau_0). 
\end{equation}
It is easy to see that $\Ueff(\tau,\tau_0)$ satisfies the equation
\begin{equation}
i \eps\frac{d}{d\tau}\Ueff(\tau,\tau_0) = \Heff(\tau) \Ueff(\tau,\tau_0), \ \Ueff(\tau_0,\tau_0) = \textrm{id},
\end{equation}
where $\Heff(\tau)$ is the \emph{effective Hamiltonian} whose form is
\begin{align}
\Heff (\tau)&=	\mathcal{U}(\tau) h_a(\tau) \mathcal{U}^*(\tau) + i \eps \dot{\mathcal{U}}(\tau) \mathcal{U}^*(\tau)
\nonumber \\
&=\begin{pmatrix} \label{Heff}
\lambda_{\alpha}(\tau) & 0\\
 0 & \lambda_{\beta}(\tau)
\end{pmatrix}
- i \eps 
\begin{pmatrix}
\langle\phi_{\alpha}(\tau), \dot{\phi}_{\alpha}(\tau)\rangle & \langle\phi_{\beta}(\tau),\dot{\phi}_{\alpha}(\tau)\rangle\\
\langle\phi_{\alpha}(\tau), \dot{\phi}_{\beta}(\tau)\rangle &\langle\phi_{\beta}(\tau), \dot{\phi}_{\beta}(\tau)\rangle
\end{pmatrix}.
\end{align}

Theorem~\ref{th: adiabatic} implies the following.
\bt \label{effective}
Assume that $\{\lambda_j,\lambda_{j+1}\}$ is a \mariospectrum\ on $\omega'$ and let $\uu : [t_0,t_f] \rightarrow \omega'$ be a $\con^2$ curve 
such that there exists a $\con^1$-varying  basis of $\Hdue(\cdot)$ made of eigenstates of $\h(\cdot)$.
Then there exists a constant $C$ such that
\[
\|\left( U^{\eps}(\tau,\tau_0)-\mathcal{U}^*(\tau) \Ueff (\tau,\tau_0) \mathcal{U}(\tau_0)\right)\Pdue(\tau_0) \|  \leq C \eps (1+|\tau-\tau_0|) 
\]
for every $\tau,\tau_0\in [t_0,t_f]$. 
\et

\subsection{Regularity of eigenstates} 
\label{s-reg}
Classical results (see \cite{reed_simon}) say that the map  $\uu\mapsto P_{\uu}$, where $P_{\uu}$ is the spectral projection relative to a \mariospectrum, is analytic on $\omega$.
In particular, eigenstates relative to simple eigenvalues can be chosen analytic with respect to $\uu$. 

Similar results hold also for intersecting eigenvalues, provided that the Hamiltonian depends on one parameter and is analytic. In particular, if $\Sigma$ is a \mariospectrum\ on $\omega$ and $\uu:I\to\omega$ is analytic, then it is possible to find 
two families of analytic functions
$\Lambda_j:I \to \mathbb{R}$
and $\Phi_j:I \to \mathcal{H} , \ j=0,\ldots,k$,    such that
for any $t$ in $I$ we have  $\Sigma(\uu(t))=\{\Lambda_0(t),\ldots,\Lambda_k(t)\}$ and 
$(\Phi_0(t),\ldots,\Phi_k(t))$ is an orthonormal basis of
corresponding eigenstates (see \cite{katino}, \cite[Theorem XII.13]{reed_simon}).

Moreover, we can easily find conditions on the derivatives of the functions $\Lambda_l,\Phi_l$: indeed, consider 
a $\con^1$ curve $\uu: I \to \mathbb{R}^2$ such that there exist
two families of $\con^1$ functions
$\Lambda_l:I \to \mathbb{R}$
and $\Phi_l:I \to \mathcal{H} , \ l=0,\ldots,k$, which 
for any $t\in I$, correspond to the eigenvalues and the (orthonormal) eigenstates of $H(\uu(t))$. 
 
By direct computations we obtain that for all $t\in I$ the following equations hold: 
\begin{gather}
\dot{\Lambda}_l(t) = \langle \Phi_l(t), \left(\dot{u}_1(t) H_1+\dot{u}_2(t) H_2 \right) \Phi_l(t) \rangle \label{obs1}
\\
(\Lambda_m(t)-\Lambda_l(t)) \:\langle \Phi_l(t),\dot{\Phi}_m(t)\rangle=
\langle \Phi_l(t), \left(\dot{u}_1(t) H_1+\dot{u}_2(t) H_2 \right) \Phi_m(t) \rangle.\label{obs2}
\end{gather}

An immediate consequence of \eqref{obs1} is that the eigenvalues $\lambda_l$ are Lipschitz with respect to $t$.

Let  $\bar \uu$ be a conical intersection between $\lam_j(\uu)$ and $\lam_{j+1}(\uu)$.
Consider the straight line $r_{\theta}(t)=\bar \uu +t(\cos \theta,\sin \theta), \ t\geq 0$.
Then \eqref{obs2} implies that
\beq \label{limitzero}
\lim_{t\rightarrow 0^+} \langle \phi_j(r_{\theta}(t)), (\cos \theta H_1+\sin \theta H_2) \phi_{j+1}(r_{\theta}(t)) \rangle=0.
\eeq

\section{Conical intersections} \label{sconical}

From now on, we assume that the Hamiltonian satisfies hypothesis {\bf (H1)}.
Following Remark 
\ref{reale}, we always choose the eigenfunctions of $H(\uu)$ whose components are real with respect to the basis $\{\chi_l\}_l$ defined in hypothesis {\bf (H1)}. In particular, this ensures that the values $\langle \phi_l(\uu), H_0 \phi_m(\uu) \rangle,\ \langle \phi_l(\uu), H_1 \phi_m(\uu) \rangle$ and $\langle \phi_l(\uu), H_2 \phi_m(\uu) \rangle,\ l,m=0,\ldots,k$, are real for any $\uu$. 

In this section, we investigate the features of conical intersections and provide also a criterion for checking if an intersection between two eigenvalues is conical. 
First of all we notice that Definition~\ref{conical} can be reformulated by saying that 
an intersection $\bar \uu$ between the eigenvalues $\lam_{j}$ and $\lam_{j+1}$ is conical if and only if there exists $c>0$ such that, for every straight line $r(t)$
with $r(0)=\bar\uu$, it holds
\[\frac{d}{dt} \Big|_{t=0^+}\Big[\lam_{j+1}(r(t))-\lam_{j}(r(t)) \Big]\geq c.\]

Moreover, the following result guarantees that \eqref{formcono} holds true in a neighborhood of a conical intersection. It follows directly from the Lipschitz continuity of the eigenvalues.

\bl
\label{ovvio}
Let $\bar \uu$ a conical intersection between $\lam_j$ and $\lam_{j+1}$. Then there exists a suitably small neighborhood $U$ of $\bar\uu$ and $C>0$ such that  
\bqn
\lam_{j+1}(\uu)-\lam_j(\uu)\geq C| \uu - \bar \uu|,\ \forall \uu\in U.
\label{tesi}
\eqn
\el

Let us now introduce the following matrix, which plays a crucial role in our controllability result.	

\bdeff
Let $\psi_1,\psi_2$ be a pair of elements of $\H$.
The {\it \mariomatrix\ associated with $(\psi_1,\psi_2)$} is
\beq
{\cal M}(\psi_1,\psi_2)=\left(\ba{cc} \langle\psi_1,H_1\psi_2\rangle & \frac12 \big(\langle\psi_2,H_1\psi_2\rangle - \langle\psi_1,H_1\psi_1\rangle\big) \\
\langle\psi_1,H_2\psi_2\rangle & \frac12\big(\langle\psi_2,H_2\psi_2\rangle - \langle\psi_1,H_2\psi_1\rangle\big) \ea\right).
\eeq
\edeff

\bl
The function $(\psi_1,\psi_2) \mapsto |\det{\cal M}(\psi_1,\psi_2)|$ is invariant under orthogonal transformation of the argument, that is if $( \widehat{\psi}_1, \widehat{\psi}_2)^T= \mathcal{O}(\psi_1, \psi_2)^T$ for a pair $\psi_1,\psi_2$ of orthonormal elements of $\H$ and $\mathcal{O}\in \textsf{O(2)} $, 
then 
one has $|\det \mathcal{M}(\widehat{\psi}_1,\widehat{\psi}_2)|=|\det \mathcal{M}(\psi_1,\psi_2)|$.
\label{tuttugulu}
\el
\proof
We set $\mathcal{O}=\left(\begin{smallmatrix} \cos\al &  \sin\al \\ -\varsigma\sin\al & \varsigma \cos\al\end{smallmatrix}\right)$, where $\varsigma=\pm 1$. A direct computation shows that \[\mathcal{M}(\widehat{\psi}_1,\widehat{\psi}_2)=\mathcal{M}(\psi_1,\psi_2)\left(\ba{cc} \cos2\al &- \sin2\al \\ \sin2\al & \phantom{-}\cos2\al\ea\right) 
\begin{pmatrix}
\varsigma  & 0\\
0 & 1
\end{pmatrix}
,\]
which immediately leads to the thesis. 
\eproof

The following result characterizes conical intersections in terms of the \mariomatrix.

\bp 
\label{iff conical}
Assume that $\{\lam_j, \lam_{j+1}\}$ is a \mariospectrum, and $\lam_j(\bar \uu) = \lam_{j+1}(\bar \uu)$.
Let $\{\psi_1,\psi_2\}$ be an orthonormal basis of the eigenspace associated with the double eigenvalue. Then $\bar\uu$ is a conical intersection if and only if ${\cal M}(\psi_1,\psi_2)$ is nonsingular.
\ep

\noindent 
\proof
Let $r_{\theta}(t)=\bar\uu+t(\cos\theta,\sin\theta)$ 
and let $\phi_j^{\theta},\phi_{j+1}^{\theta}$ be the limits of $\phi_j(r_{\theta}(t)),\phi_{j+1}(r_{\theta}(t))$ as $t\rightarrow 0^+$ (recall that the eigenfunctions $\phi_j,\phi_{j+1}$ can be chosen analytic along $r_{\theta}$ for $t\geq 0$). 
Assume that for any $\eps>0$ there exists $\theta_{\eps}$ such that
\[\frac{d}{dt} \Big|_{t=0^+} \Big[\lam_{j+1}(r_{\theta_{\eps}}(t))-\lam_{j}(r_{\theta_{\eps}}(t)) \Big] \leq \eps,\] 
that is, by \eqref{obs1}, $\cos \theta_{\eps} \left( \langle \phi_j^{\theta_{\eps}}, H_1\phi_j^{\theta_{\eps}} \rangle - \langle \phi_{j+1}^{\theta_{\eps}}, H_1 \phi_{j+1}^{\theta_{\eps}} \rangle \right) +\sin \theta_{\eps} \left( \langle \phi_j^{\theta_{\eps}}, H_2  \phi_j^{\theta_{\eps}} \rangle - \langle \phi_{j+1}^{\theta_{\eps}}, H_2 \phi_{j+1}^{\theta_{\eps}} \rangle \right) \leq \eps$. Moreover, by \eqref{limitzero}, we have that 
$\cos \theta_{\eps} \langle \phi_j^{\theta_{\eps}}, H_1\phi_{j+1}^{\theta_{\eps}} \rangle +\sin \theta_{\eps} \langle \phi_j^{\theta_{\eps}}, H_2\phi_{j+1}^{\theta_{\eps}} \rangle =0$.
Since 
\[\Big|\det \mathcal{M}(\phi_j^{\theta_{\eps}},\phi_{j+1}^{\theta_{\eps}})\Big|=\Big|\det \left[
\begin{pmatrix}
\cos \theta_{\eps} & \sin \theta_{\eps}\\
-\sin \theta_{\eps} & \cos \theta_{\eps}
\end{pmatrix}
\mathcal{M}(\phi_j^{\theta_{\eps}},\phi_{j+1}^{\theta_{\eps}})\right] \Big| \leq 2\eps(\|H_1\|+\|H_2\|),
\]
then, by Lemma \ref{tuttugulu} and the arbitrariness of $\eps$, we get that ${\cal M}(\psi_1,\psi_2)$ is singular. 
Thus $\bar\uu$ is a conical intersection when ${\cal M}(\psi_1,\psi_2)$ is nonsingular.

Let us now prove the converse statement: assume that $\bar \uu$ is a conical intersection and, by contradiction, that $\mathcal{M}(\phi_j^{\theta},\phi_{j+1}^{\theta})$ is singular,
where  $\phi_j^{\theta},\phi_{j+1}^{\theta}$ are defined as above. 
By definition of conical intersection, we have
\begin{equation}
\cos \beta \left( \langle \phi_j^{\beta}, H_1\phi_j^{\beta} \rangle - \langle \phi_{j+1}^{\beta}, H_1 \phi_{j+1}^{\beta} \rangle \right) +\sin \beta \left( \langle \phi_j^{\beta}, H_2  \phi_j^{\beta} \rangle - \langle \phi_{j+1}^{\beta}, H_2 \phi_{j+1}^{\beta} \rangle \right) \neq 0, \label{conint}
\end{equation}
for every $\beta \in \R$.

By  \eqref{limitzero} and \eqref{conint} with $\beta=\theta$, it turns out that the two columns of the matrix $\mathcal{M}(\phi_j^{\theta},\phi_{j+1}^{\theta})$ are not proportional. Thus $\mathcal{M}(\phi_j^{\theta},\phi_{j+1}^{\theta})$ can be singular only if its first column is null. 

For any angle $\beta$,  there exists an orthonormal matrix $\mathcal{O}=\left(\begin{smallmatrix} \cos\al & \sin\al \\ -\sin\al & \cos\al\end{smallmatrix}\right)$ such that $(\phi_j^{\beta},\phi_{j+1}^{\beta})^T=\mathcal{O}(\phi_j^{\theta},\phi_{j+1}^{\theta})^T$ and, calling $W=\cos \beta H_1+\sin \beta H_2$, we have (by \eqref{limitzero})
\begin{align*}
0=\langle \phi_{j}^{\beta}, W \phi_{j+1}^{\beta} \rangle &= (\cos \al^2 - \sin \al^2) \langle \phi_j^{\theta}, W \phi_{j+1}^{\theta} \rangle + \sin \al \cos \al \left(\langle \phi_{j+1}^{\theta}, W \phi_{j+1}^{\theta} \rangle-\langle \phi_{j}^{\theta}, W \phi_{j}^{\theta} \rangle \right) =\\
&=  \sin \al \cos \al \left(\langle \phi_{j+1}^{\theta}, W \phi_{j+1}^{\theta} \rangle-\langle \phi_{j}^{\theta}, W \phi_{j}^{\theta} \rangle \right).
\end{align*}
If $\langle \phi_{j+1}^{\theta}, W \phi_{j+1}^{\theta} \rangle-\langle \phi_{j}^{\theta}, W \phi_{j}^{\theta} \rangle = 0$, the matrix $(\langle \phi_l^{\theta}, W \phi_m^{\theta}\rangle )_{l,m=j,j+1}$ is diagonal and proportional to the identity. Hence the same is true for  $(\langle \phi_l^{\beta}, W \phi_m^{\beta}\rangle )_{l,m=j,j+1}$.  This contradicts \eqref{conint}, so that it must be $\sin \alpha \cos \alpha =0$, that is, the limit basis is unique and therefore it must be equal to $\{\phi_j^{\theta},\phi_{j+1}^{\theta}\}$ (up to phases).

Let us now consider the straight line $r_{\beta}$ with
\[
\tan \beta=\frac{\langle \phi_{j+1}^{\theta}, H_1 \phi_{j+1}^{\theta} \rangle-\langle \phi_{j}^{\theta}, H_1 \phi_{j}^{\theta} \rangle}{\langle \phi_{j+1}^{\theta}, H_2 \phi_{j+1}^{\theta} \rangle-\langle \phi_{j}^{\theta}, H_2 \phi_{j}^{\theta} \rangle}.
\]
Since, as proved above, the limit basis along $r_{\beta}$ is $\{\phi_j^{\theta},\phi_{j+1}^{\theta}\}$
we have that
$\langle \phi_j^{\theta} , (\cos \beta H_1 + \sin \beta H_2)\phi_j^{\theta}  \rangle = \langle \phi_{j+1}^{\theta} , (\cos \beta H_1 + \sin \beta H_2) \phi_{j+1}^{\theta} \rangle$.
By \eqref{obs1}, this contradicts \eqref{conint}, proving that $\mathcal{M}(\phi_j^{\theta},\phi_{j+1}^{\theta})$ is nonsingular.
\eproof
%

As noticed above, for any analytic curve that reaches a conical intersection it is possible to choose continuously the  eigenstates along the curve.
A peculiarity  of conical intersections is that, when approaching the singularity from different directions, the eigenstates corresponding to the intersecting eigenvalues have different limits, and the dependence of such limits from the direction can be explicitly computed, as shown in the following result.

\newcommand{\bch}{\boldsymbol{\chi}}

\bp
\label{tetha}
Let $\bar \uu$ be a conical intersection between the eigenvalues $\lam_j,\lam_{j+1}$ and let $\phi_j^0,\phi_{j+1}^0$ be the limits as $t \to 0^+$ of the eigenstates $\phi_j(r_0(t) ),\phi_{j+1}(r_0(t) )$, for $r_0(t)=\bar \uu+ (t,0)$. Consider, for any $\al\in [0,2\pi)$, the straight line $r_\al (t)=\bar \uu+(t\cos\al,t\sin\al)$. Then, up to a sign,
the eigenstates $\phi_j(r_{\al}(t) ),\phi_{j+1}(r_{\al}(t) )$ have limits
\begin{align}
\phi_j^\al=& \cos\bth(\al)\phi_j^0+\sin\bth(\al)\phi_{j+1}^0  \label{transangolo1}\\
\phi_{j+1}^\al=&-\sin\bth(\al)\phi_j^0+\cos\bth(\al)\phi_{j+1}^0, \label{transangolo2}
\end{align}
 where $\bth$ is a monotone $\con^1$ function defined on $[0,2\pi)$ with $\bth(0)=0$.
Depending on the initial choice of  $\phi_j^0,\phi_{j+1}^0$ the range of $\bth$ is either $[0,\pi)$ or $(-\pi,0]$.
Moreover, $\bth(\cdot)$ satisfies the equation
\beq \label{thetaalpha}
\begin{pmatrix}
\cos \alpha, \sin \alpha
\end{pmatrix} 
\mathcal{M}(\phi_j^0,\phi_{j+1}^0)
\begin{pmatrix}
\cos 2\bth(\al) \\
\sin 2\bth(\al)
\end{pmatrix} =0.
\eeq
\ep

\proof
Let us write $\phi_j^\al,\phi_{j+1}^\al$ as in \eqref{transangolo1}-\eqref{transangolo2}. 
Then $\bth(\alpha)$ satisfies
\begin{align}
0&=\langle \phi_j^\al,(\cos \alpha H_1+\sin \alpha H_2)\phi_{j+1}^\al \rangle \nonumber \\
&=  \cos 2\bth(\alpha )\langle \phi_j^0,(\cos \alpha H_1+\sin \alpha H_2)\phi_{j+1}^0 \rangle + \nonumber  \\
&+\frac{1}{2} \sin 2\bth(\alpha) \left(\langle \phi_{j+1}^0,(\cos \alpha H_1+\sin \alpha H_2)\phi_{j+1}^0 \rangle-\langle \phi_j^0,(\cos \alpha H_1+\sin \alpha H_2)\phi_{j}^0 \rangle \right)  \nonumber \\
&=  \begin{pmatrix}
\cos \alpha, \sin \alpha
\end{pmatrix} 
\mathcal{M}(\phi_j^0,\phi_{j+1}^0)
\begin{pmatrix}
\cos 2\bth(\al) \\
\sin 2\bth(\al)
\end{pmatrix},  \nonumber \label{limitzero=0}
\end{align}
%
proving \eqref{thetaalpha}.
Equation \eqref{thetaalpha} has exactly four solutions for any value of $\alpha$, differing one from the other by  multiples of $\pi/2$. 
By the Implicit Function Theorem, it turns out that each of them is a $\con^1$ monotone function defined on $[0,2\pi)$.

We define $\bth(\cdot)$ as the one that satisfies $\bth(0)=0$. We are left to prove that the range of $\bth$ is $[0,\pi)$ or $(-\pi,0]$.
We first observe that when $\alpha=\pi$ the possible solutions of equation \eqref{thetaalpha} are multiples of $\pi/2$. If $|\bth(\pi)| > \pi/2$, then by continuity there should exist $\bar \alpha \in (0,\pi)$ 
with $|\bth(\bar \alpha)|=\pi/2$. This is impossible because of equation \eqref{thetaalpha}. 
Thus $\bth$ maps $[0,\pi]$ into $[0,\pi/2]$ or $[-\pi/2,0]$ and, by symmetry, the claim is proved. 
\hfill $\Box$

\brem
From Proposition \ref{tetha} it is straightforward to see that it is not possible to define continuously the eigenstates $\phi_j,\phi_{j+1}$ of $H(\uu)$ on a closed path that encloses the singularity: after a complete turn, a change of sign appears. 
\erem

\section{\curvepao}
\label{una-a-caso}
Throughout this section we will assume that $\{\lam_j,\lam_{j+1}\}$ is a \mariospectrum\ on some open domain $\omega$, and that $0 \in \omega$ is a conical intersection between the eigenvalues. Without loss of generality, in the following we always assume that $0$ is the only intersection between $\lam_j$, $\lam_{j+1}$ in $\omega$.

Following Section \ref{subs ad}, the effective Hamiltonian $\Heff$, defined as in \eqref{Heff}, (approximately) describes the dynamics inside the eigenspaces associated with $\lam_j,\lam_{j+1}$, for $\uu$ slowly varying inside $\omega$. 

When integrating the effective Hamiltonian, the second term in \eqref{Heff} gives a total contribution that a priori is of order $O(1)$. In particular the contribution of the non-diagonal terms of $\Heff$ induce a (a priori) non-negligible probability transfer between the two levels.


To tackle this issue we consider trajectories satisfying the following dynamical system
\beq
\left\{\ba{l} 
\dot u_1 =  -\langle\phi_{j},H_2\phi_{j+1}\rangle\\ 
\dot u_2 =  \langle\phi_{j}, H_1\phi_{j+1}\rangle.
\ea\right.
\label{pao}
\eeq
Notice that the right-hand side of \eqref{pao} can be taken real-valued under hypothesis {\bf (H1)}. 
It is defined up to a sign, because of the freedom in the choice of the sign of the eigenstates.
Nevertheless, { locally around points where $\lam_j\neq\lam_{j+1}$,} it is possible to choose the sign in such a way that the right-hand side of \eqref{pao}  is smooth, and,
from equation \eqref{obs2}, we see that $\langle \phi_j(\gamma(t)),\dot{\phi}_{j+1}(\gamma(t))\rangle=0$ along any integral curve $\gamma$ of \eqref{pao}. Here and in the following we use the notation $\dot\phi (\gamma(\cdot))$ to denote $\frac{d}{dt}(\phi(\gamma(\cdot)))$.

Let now $\mathcal{H}^{\mathbb{R}}$ be the real Hilbert space generated by the basis $\{\chi_j\}_j$ defined in Remark~\ref{reale}, and let $ \mathrm{Gr}_2({\cal H}^{\mathbb{R}})$ be the 2-Grassmannian of ${\cal H^{\mathbb{R}}}$, i.e. the set of all two-dimensional subspaces of ${\cal H^{\mathbb{R}}}$. This set has a natural structure of metric space defined by the distance $d(W_1,W_2)=\|P_{W_1}-P_{W_2}\|$, where  $P_{W_1},P_{W_2}$ are the orthogonal projections on the two-dimensional subspaces $W_1,W_2$.  Lemma~\ref{tuttugulu} allows us to define the function 
\begin{align}
\hat F : \mathrm{Gr}_2({\cal H}^{\mathbb{R}}) &\rightarrow \R \label{FF} \\
W &\mapsto |\det {\cal M}(v_1,v_2)| \nonumber,
\end{align}
 where $\{v_1,v_2\}$ is any orthonormal basis of $W \in \mathrm{Gr}_2({\cal H^{\mathbb{R}}})$.
It is straightforward to see that  $\hat F$ is continuous.
	

Let $P_{\uu}$ be the spectral projection associated with  the pair $\{\lam_j(\uu),\lam_{j+1}(\uu)\}$. 
We know from Section~\ref{s-reg}  that $P_{\uu}$ is analytic on $\omega$. Therefore $\uu \mapsto P_{\uu}\H \cap \H^{\mathbb{R}}$ is continuous in Gr$_2(\H^{\mathbb{R}})$.
Let now $F(\uu):=|\det {\cal M}(\phi_j(\uu),\phi_{j+1}(\uu))|$. Since $F(\uu)=\hat F( P_{\uu}\H \cap \H^{\mathbb{R}})$  and by Proposition~\ref{iff conical}  we get the following result.

\bl
The function $\uu \mapsto F(\uu)$ is well defined and continuous in $\omega$. In particular $F$ is different from $0$ in a neighborhood of $\uu=0$. 
\label{det-cont}
\el

Without loss of generality, we assume from now on that $F$ is different from zero on $\omega$. 

\bl
There exists a $\con^{\infty}$ choice of the right-hand side of \eqref{pao} in $\omega\setminus\{0\}$ such that, if $\uu(\cdot)$ is a corresponding solution, then 
\begin{equation} \label{ruby}
\frac{d}{dt}\Big[\lam_{j+1}(\uu(t))-\lam_{j}(\uu(t)) \Big]=-F(\uu(t))
\end{equation}
on $\omega\setminus\{0\}$.
\label{diff}
\el
\proof Observe that 
\[
\frac{d}{dt}\Big[\lam_{j+1}(\uu(t))-\lam_{j}(\uu(t)) \Big]=
\dot u_1 \big(\langle\phi_{j+1},H_1\phi_{j+1}\rangle- \langle\phi_j,H_1\phi_j\rangle\big)
+\dot u_2 \big(\langle\phi_{j+1},H_2\phi_{j+1}\rangle-\langle\phi_j,H_2\phi_j\rangle\big).
\]
This expression, evaluated along the solutions of \eqref{pao}, is equal either to $F(\uu(t))$ or to $-F(\uu(t))$, depending on the choice of the sign in \eqref{pao}. Since $F(\uu)\neq 0$ on $\omega$, there exists a unique choice of this sign such that equation \eqref{ruby} is satisfied. The local smoothness of the eigenfunctions ensures that this choice is smooth.  
\eproof

We now define the \emph{\campopao}, 
denoted by $\Xpao$, as
the smooth vector field on $\omega\setminus\{0\}$ identified by the preceding lemma. 
Its integral curves are $\Cspet$ in $\omega\setminus\{0\}$.
Moreover,
its norm is equal to the norm of the first row of $\mathcal{M}(\phi_j,\phi_{j+1})$, and therefore bounded both from above and from below by positive constants in $\omega\setminus\{0\}$.

By considering $\lambda_{j+1}(\uu)-\lambda_{j}(\uu)$ as a local Lyapunov function, the above results lead to the following proposition. 
\bp \label{arrivanoalcentro}
There exists a punctured neighborhood $U$ of $0$  such that all the integral curves of $\Xpao$ starting from $U$ reach the origin in finite time.
\ep

Our purpose now is to prove that each of these curves admits a $\con^\infty$ extension up to the singularity.
As a preliminary result we get the following.

\bp \label{limitati}
Let  $\bar\uu=0$ be a conical intersection with $\lam_j(0)=\lam_{j+1}(0)$, and let the map $\uu \mapsto \{\lam_j(\uu),\lam_{j+1}(\uu)\}$ be a \mariospectrum\ on a neighborhood of $0$. Then, for any $C_1>0$, there exist a neighborhood $I$ of $t=0$ and $C_2>0$ such that for any $\con^2$ trajectory $\gamma(\cdot)$ with $\gamma(0)=0$, $|\dot\gamma(0)|=1$ and $\|\ddot{\gamma}\|_{L^{\infty}(I)}\leq C_1$, one has $\|\dot{\phi}_l(\gamma(t))\|\leq C_2,\ l =j,j+1$, for every $t\in I\setminus\{0\}$.
\ep
\proof 
Let us assume without loss of generality that $\lam_j(0)=\lam_{j+1}(0)=0$.
For $t\neq 0$ define $\rho(t)=\dot\gamma(t)-\gamma(t)/t$. 
Notice that $\rho(t)=\frac12 \ddot \gamma(0)t+o(t)$.

By \eqref{obs2}, we have 
\beq
\langle\phi_j(\gamma(t)),\dot{\phi}_{j+1}(\gamma(t))\rangle=\frac{\langle\phi_j,(\gamma_1 H_1+  \gamma_2 H_2)\phi_{j+1}\rangle}{t(\lambda_{j+1}-\lambda_{j})}+\frac{\langle\phi_j,(\rho_1 H_1+ \rho_2 H_2)\phi_{j+1}\rangle}{\lambda_{j+1}-\lambda_{j}}.
\label{2term}
\eeq
Notice that ${\langle\phi_j(\uu),(u_1 H_1+ u_2 H_2)\phi_{j+1}(\uu)\rangle}=-{\langle\phi_j(\uu),H_0\phi_{j+1}(\uu)\rangle}= -\langle \phi_j(\uu)-P_0\phi_j(\uu),H_0 \big(\phi_{j+1}(\uu) - P_0\phi_{j+1}(\uu)\big)\rangle$.
Since 
\begin{align*}
\|H_0\big(\phi_{j}( \uu ) - P_0\phi_{j}(\uu)\big)\| &=  \| \lam_{j}(\uu)\phi_{j}(\uu)   - u_1 H_1\phi_{j}(\uu) - u_2 H_2\phi_{j}(\uu)\|\nonumber\\
&\leq  | \lam_{j}(\uu)|  + |\uu|(\|H_1\|+\|H_2\|)\nonumber\\
&\leq  \sup_{\boldsymbol{v} \in \omega} |\la\phi_{j}(\boldsymbol{v}),(u_1H_1+u_2H_2)\phi_{j}(\boldsymbol{v})\ra| +  |\uu|(\|H_1\|+\|H_2\|)\nonumber\\
&\leq  2(\|H_1\|+\|H_2\|) |\uu| \nn 
\end{align*}
and by smoothness of the projector, 
we get  that $\big|\langle \phi_j(\uu)-P_0\phi_j(\uu),H_0 \big(\phi_{j+1}(\uu) - P_0\phi_{j+1}(\uu)\big)\rangle\big|\leq 2C(\|H_1\|+\|H_2\|)|\uu|^2$, for a suitable $C>0$. Being $|\gamma(t)|=O(t)$ and $\lambda_{j+1}(\uu)-\lambda_{j}(\uu)>c|\uu|$ (Lemma~\ref{ovvio}), we deduce that the  modulus of the first term in the right-hand side of \eqref{2term} is uniformly bounded. 
The uniform bound of the second term is a trivial consequence of the fact that $|\rho(t)|=O(t)$ and that $|\gamma(t)|\geq \bar{c} |t|$, for some $\bar{c}>0$, if $t$ is small enough. Thus $|\langle\phi_j(\gamma(t)),\dot\phi_{j+1}(\gamma(t))\rangle|$ is uniformly bounded.

Let us write  $P_{\uu}^{\bot}=\mathrm{id}-P_{\uu}$. Since $P_{\uu}^{\bot}$ commutes with $H(\uu)$, one has
\[(H(\gamma(t))-\lam_{j+1}(\gamma(t))\mathrm{id})P_{\gamma(t)}^{\bot}\dot\phi_{j+1}(\gamma(t))=-P_{\gamma(t)}^{\bot} (\dot\gamma_1(t) H_1+  \dot\gamma_2(t) H_2)\phi_{j+1}(\gamma(t)).\]
Since $H(\uu)-\lam_{j+1}(\uu)\mathrm{id}$ is invertible on $P_{\uu}^{\bot}\H$ with uniformly bounded inverse on $\omega$, we get that $\|P_{\gamma(t)}^{\bot}\dot\phi_{j+1}\|$ is uniformly bounded on $I\setminus\{0\}$. Thus we obtain that $\|\dot\phi_{j+1}\|$ is bounded, uniformly on the set of curves $\gamma(\cdot)$ satisfying the assumptions of the proposition.
The same holds for $\|\dot\phi_{j}\|$.
\eproof

\bcc \label{cursore}
Let  $\bar\uu=0$ be a conical intersection with $\lam_j(0)=\lam_{j+1}(0)$, and let the map $\uu \mapsto \{\lam_j(\uu),\lam_{j+1}(\uu)\}$ be a \mariospectrum\ on a neighborhood of $0$. Denote by $\tphi_l(\rho,\theta)$ the eigenstate $\phi_l(\rho\cos \theta,\rho \sin \theta), \ l=j,j+1$, where $(\rho,\theta)$ are angular coordinates  around $0$, i.e. $\rho=|u|$ and $\theta=\arctan \frac{u_2}{u_1}$. Set $\tphi_l(0,\theta)=\lim_{\rho\rightarrow 0^+}\tphi_l(\rho,\theta)$.
Then the function $\tphi_l(\rho,\theta)$  is continuous in $[0,R]\times [0,2\pi]$, for some $R>0$ and $l=j,j+1$.    
\ecc

\proof If $\rho>0$, the function $\tphi_l(\rho,\theta)$ can be defined  continuously. Moreover, the function $\theta \mapsto \tphi_l(0,\theta)$ is uniformly continuous, thanks to Proposition \ref{tetha}. 

Let us now consider a sequence $(\rho_k,\theta_k)$ converging to $(0,\bar\theta)$. Then we have
\begin{align*}
|\tphi_l(\rho_k,\theta_k)-\tphi_l(0,\bar \theta)| &\leq |\tphi_l(\rho_k,\theta_k)-\tphi_l(0,\theta_k)|+|\tphi_l(0,\theta_k)-\tphi_l(0,\bar \theta)|\\
&\leq C_1 \rho_k+|\tphi_l(0,\theta_k)-\tphi_l(0,\bar \theta)|,
\end{align*}
where $C_1$ comes from Proposition~\ref{limitati} and the second term goes to zero as $k$ goes to infinity.  
\eproof

\bp
The eigenstates $\phi_j,\phi_{j+1}$ can be extended continuously to the singularity along the integral curves of $\Xpao$, and, in a small enough punctured neighborhood of $\uu=0$, the integral curves of $\Xpao$ admit a $\con^1$ extension up to the singularity included.
\ep

\proof
We prove that the scalar product $\Xpao \cdot (-u_2,u_1)^T/|\uu|$ goes to $0$ as $|\uu|\to 0$, that is, the tangent to the curve has limit when $\uu$ approaches zero.
This, together with Corollary \ref{cursore}, implies that the eigenstates $\phi_j,\phi_{j+1}$ are continuous along the integral curves of $\Xpao$, and then the vector field $\Xpao$ itself is continuous along its integral curves, up to the singularity included. Therefore 
its integral curves admit a $\con^1$ extension up to the singularity.

To prove that $\Xpao \cdot (-u_2,u_1)^T/|\uu|$ goes to $0$ as $|\uu|\to 0$, we show that there exists a constant $C>0$ such that
\beq
\kappa(\uu): = | \Xpao(\uu) \cdot (-u_2,u_1)  | \leq C|\uu|^2\,.\label{stima}
\eeq
Since 
$\kappa(\uu)=|\langle\phi_j(\uu),(u_1 H_1+u_2 H_2)\phi_{j+1}(\uu)\rangle|$, the thesis comes from the estimates in the proof of Proposition \ref{limitati}. \eproof

We recall that, since integral curves of the \campopao\ $\Xpao$ are $\con^1$, then the spectral projection $P_{\uu}$ associated with the pair $\{\lam_j(\uu),\lam_{j+1}(\uu)\}$ is $\con^1$ along each of them. This permits to prove the following result.


\bp\label{spettacolo}
For any integral curve $\gamma : [-\eta,0] \rightarrow \omega$
of $\Xpao$ with $\gamma(0)=0$ there exists a choice of an orthonormal basis of the eigenspace associated with the double eigenvalue $\lambda_j(\gamma(0))=\lambda_{j+1}(\gamma(0))$ that makes the eigenstates $\phi_j(\gamma(t)),\phi_{j+1}(\gamma(t))$ $\con^1$ on $[-\eta,0]$.
\ep 

\proof
We notice that on the integral curves of $\Xpao$ the eigenstates relative to the eigenvalues $\lam_j,\lam_{j+1}$ satisfy the equation $P_{\gamma(t)} \dot{\phi}_j(\gamma(t))=P_{\gamma(t)} \dot{\phi}_{j+1}(\gamma(t))=0$, which implies
\begin{equation} \label{der-spettacolo}
\dot{P}_{\gamma(t)} \phi_j(\gamma(t)) =\dot{\phi}_j(\gamma(t))  \quad
\dot{P}_{\gamma(t)} \phi_{j+1}(\gamma(t)) =\dot{\phi}_{j+1}(\gamma(t)) 
\end{equation}
for $t\in [-\eta,0)$. The thesis follows from the continuity of $\dot{P}_{\gamma(t)}, \phi_j(\gamma(t)), \phi_{j+1}(\gamma(t)) $ on $[-\eta,0]$. \hfill $\Box$

\bcc
Let $\gamma : [-\eta,0] \rightarrow \omega$ be an integral curve of $\Xpao$ with $\gamma(0)=0$. Then $\gamma(\cdot)$ and the
eigenstates $\phi_j(\gamma(\cdot)),\phi_{j+1}(\gamma(\cdot))$ defined in Proposition~\ref{spettacolo} are $\Cspet$ on $[-\eta,0]$.
\ecc

\proof
Extend $\Xpao (\gamma(t))$ by setting 
\[
\Xpao (\gamma(0))=
\begin{pmatrix} 
-\langle \phi_j(\gamma(0)), H_2\phi_{j+1}(\gamma(0))\rangle\\
\langle \phi_j(\gamma(0)),H_1\phi_{j+1}(\gamma(0)) \rangle \end{pmatrix},\] 
where $ \phi_j(\gamma(0)), \phi_{j+1}(\gamma(0))$ denote the limits of the eigenstates as defined in Proposition~\ref{spettacolo}. Then  $\Xpao(\gamma(\cdot))$ is
$\con^1$ on $[-\eta,0]$, which implies that $\gamma(\cdot)$ is $\con^2$ on $[-\eta,0]$. 
We differentiate equation \eqref{der-spettacolo} to prove that $\phi_j(\gamma(\cdot)),\phi_{j+1}(\gamma(\cdot))$ are $\con^2$ on $[-\eta,0]$.
Repeating recursively the argument we prove the thesis.  
\hfill $\square$

We stress that, thanks to Proposition~\ref{spettacolo}, if we define the adiabatic Hamiltonian $h_a(\tau)=H(\gamma(\tau)) -i \eps P_{\gamma(\tau)} \dot{P}_{\gamma(\tau)}-i \eps P^{\bot}_{\gamma(\tau)} \dot{P}^{\bot}_{\gamma(\tau)},\ \tau=\eps t$, along integral curves of $\Xpao$, then it is possible to define the associated effective Hamiltonian, as in equation~\eqref{Heff}. 

The following result is crucial to our controllability strategy.

\bp
\label{surj}
For every unit vector ${\bf w}$ in $\R^2$ there exists an integral curve $\gamma : [-\eta,0] \rightarrow \omega$ of $\Xpao$ with $\gamma(0)=0$ such that
\[
\lim_{t\rightarrow 0^-} \frac{\dot{\gamma}(t)}{\|\dot{\gamma}(t)\|} ={\bf w}.
\]
\ep

\proof
Equation \eqref{pao} rewrites as 
\begin{align}
\dot{\rho} &=\langle\tphi_j(\rho,\theta),(-\cos\theta H_2 +\sin\theta H_1)\tphi_{j+1}(\rho,\theta)\rangle  	\label{dotro} \\
\dot\theta &=\frac{1}{\rho}\langle\tphi_j(\rho,\theta),(\cos\theta H_1 +\sin\theta H_2)\tphi_{j+1}(\rho,\theta)\rangle.
\label{ripa}
\end{align}

On a neighborhood $U\subset \omega$ of the singularity, there exist two constants $0<c_1<c_2$ such that $c_1<|\dot{\rho}|<c_2$, and the right-hand side of \eqref{ripa} is bounded from above, by \eqref{stima}. 
We choose the sign of the functions $\tphi_j,\tphi_{j+1}$ in such a way that $\dot{\rho}<0$.

Fix $\bar \theta \in [0,2\pi]$ such that ${\bf w}=(\cos\bar\theta,\sin\bar\theta)$. Consider, for $k$ large enough, the solutions $(\rho_k(\cdot),\theta_k(\cdot))$ of \eqref{dotro}-\eqref{ripa} with
$\rho_k(0)=1/k$ and $\theta_k(0)=\bar{\theta}$, for $t$ belonging to some common interval $[-\eta,0]$, where $\eta>0$ is small enough, in order to guarantee that the solutions do not exit from $U$.
By Ascoli-Arzel\`a Theorem, up to subsequences, $(\rho_k(\cdot),\theta_k(\cdot))$ converges uniformly on $[-\eta,0]$ to some
$(\tro(\cdot),\tet(\cdot))$. 

In particular, for any $\tau \in [-\eta,0]$, $(\rho_k(\tau),\theta_k(\tau))$ converges in $U$.  
By the uniform boundedness of $\dot{\rho}$, the range of $(\rho_k(\cdot),\theta_k(\cdot))$ on $[-\eta,\tau]$ is contained in a compact subset $K\subset U\setminus \{0\}$ for every $k$. Since  the vector field is smooth on $K$, the curves $(\rho_k(\cdot),\theta_k(\cdot))$ converge uniformly on $[-\eta,\tau]$ to the solution of \eqref{dotro}-\eqref{ripa} with initial condition $\rho(\tau)=\tro(\tau)$  and $\theta(\tau)=\tet(\tau)$. Therefore for $t\in[-\eta,\tau]$ $(\tro(\cdot),\tet(\cdot))$ is a solution of \eqref{dotro}-\eqref{ripa}. 
Since $\tau$ is arbitrary, and  $\tet(0)=\lim_k \theta_k(0)=\bar\theta$, $\tro(0)=\lim_k \rho_k(0)=0$, we get the thesis.
\eproof 

We conclude this section by proving a result of structural stability of conical intersections.

\bt \label{strutto}
Assume that $H(\uu)=H_0+u_1H_1+u_2H_2$ satisfies {\bf (H0)}-{\bf (H1)} and let $\bar\uu$ be a  conical intersection for $H(\uu)$ between the eigenvalues $\lambda_j$ and $\lambda_{j+1}$. Assume moreover that $\uu \mapsto \{\lambda_j(\uu),\lambda_{j+1}(\uu)\}$ is a \mariospectrum\ in a neighborhood of $\bar\uu$.
Then for every $\eps>0$ there exists $\delta>0$ such that, if $\hat H(\uu)=\hat H_0+u_1 \hat H_1+u_2\hat H_2$ satisfies  {\bf (H0)}-{\bf (H1)} and
\beq
\|\hat H_0-H_0\|+\|\hat H_1-H_1\|+\|\hat H_2- H_2\| \leq\delta,\label{varH}
\eeq
then the operator $\hat H(\uu)$ admits a conical intersection of eigenvalues at $\hat \uu$, with $|\bar\uu -\hat \uu|\leq \eps$. 
\et
\proof Continuous dependence of the eigenvalues with respect to perturbations of the Hamiltonian ensures that, if $\delta$ is small, then $\hat H$ admits two eigenvalues $\hat\lam_j,\hat\lam_{j+1}$ close 
to $\lam_j,\lam_{j+1}$.
Moreover $\{\hat\lam_j,\hat\lam_{j+1}\}$ is separated from the rest of the spectrum, locally around $\bar\uu$.
Fix now $\eps>0$ in such a way that the vector field $\Xpao$ points inside the ball $B(\bar\uu,\eps)$ at every point of its boundary (this is possible because of \eqref{stima}) and $F(\uu) \geq c>0$ on $B(\bar\uu,\eps)$.
If $\delta$ is small enough then $\hat\lam_j\neq\hat\lam_{j+1}$ on $\partial B(\bar\uu,\eps)$.
Similarly, since the conicity matrix $\mathcal{M}$ varies continuously with respect to $H_1,H_2$, and by continuity of the function $\hat F$ defined in \eqref{FF}, we can take $\delta$ small enough such that $|\det \mathcal{M}|\geq c/2$  for any perturbed Hamiltonian.
This allows us to define, whenever $\hat\lam_j\neq\hat\lam_{j+1}$,  the non-mixing field $\hat\Xpao$ associated with $\hat H$ and corresponding to the band  $\{\hat\lam_j,\hat\lam_{j+1}\}$; as in Lemma~\ref{diff}, we choose $\hat\Xpao$ in such a way that  
the time derivative of $\hat\lam_{j+1}-\hat\lam_{j}$ along the integral curves of  $\hat\Xpao$ is smaller than $-c/2$ and $\hat\Xpao$ is smooth.
In addition, by the uniform continuity on $\partial B(\bar\uu,\eps)$ of the eigenfunctions with respect to perturbations of the Hamiltonian, if $\delta$ is small enough, then $\hat\Xpao$  points inside $B(\bar\uu,\eps)$ at every point of $\partial B(\bar\uu,\eps)$.  

Fix an Hamiltonian $\hat H(\cdot)$ satisfying \eqref{varH}.  Any trajectory $\hat\gamma(\cdot)$ of $\hat\Xpao$ starting from $B(\bar\uu,\eps)$ 
 remains inside $B(\bar\uu,\eps)$ in its interval of definition and reaches in final time a point $\hat \uu$ corresponding to a double eigenvalue  $\hat\lam_j(\hat\uu) = \hat\lam_{j+1}(\hat\uu)$. 
 The conclusion follows from Proposition \ref{iff conical}.
\eproof

\section{Proof of Theorem \ref{asc-paolo}} \label{prova}
Based on Proposition~\ref{surj}, we consider below trajectories of the following kind: given a conical singularity $\uu$ and a pair of unit vectors ${\bf w}_1, {\bf w}_2 \in\R^2$, we concatenate the integral curve of  $\Xpao$ arriving at $\uu$ with direction ${\bf w}_1$ and the integral curve of  $-\Xpao$ exiting $\uu$ with direction ${\bf w}_2$.
Even if we are not using this fact in the paper, it turns out that, if ${\bf w}_1={\bf w}_2$, then such curve is globally $\Cspet$.

\bp \label{prop-2level}
Let $\uu=0$ be a conical intersection between the eigenvalues $\lam_j,\lam_{j+1}$ and let $\phi_j^0,\phi_{j+1}^0$ be limits as $\tau \to 0^+$ of the eigenstates $\phi_j(r(\tau) ),\phi_{j+1}(r(\tau) )$, respectively, for $r(\tau)=(\tau,0)$. Let $\gamma:[0,1]\rightarrow \omega$ be a piecewise $\con^{\infty}$ curve such that $\gamma(\tau_0)=0$ for some $\tau_0\in (0,1)$, $\dot{\gamma}(\tau)=\Xpao(\gamma(\tau))$ in $[0,\tau_0]$ and $\dot{\gamma}(\tau)=-\Xpao(\gamma(\tau))$ in $[\tau_0,1]$. Define  $\al_-$, $\al_+$ by
\beq
\lim_{\tau\to\tau_0^-}\frac{\dot\gamma(\tau)}{\|\dot\gamma(\tau)\|} = -(\cos\al_-,\sin\al_-)\,,\quad 
\lim_{\tau\to\tau_0^+}\frac{\dot\gamma(\tau)}{\|\dot\gamma(\tau)\|} = (\cos\al_+,\sin\al_+).
\label{eq-2level}
\eeq
Then there exists $C>0$ such that, for any  $\varepsilon >0$,
\bqn 
\|\psi(1/\eps)- p_1e^{i\vartheta_j}\phi_j(\gamma(0))-p_2e^{i\vartheta_{j+1}}\phi_{j+1}(\gamma(0))\|\leq C \eps
\eqn 
where $\vartheta_j,\ \vartheta_{j+1}\in \mathbb{R}$, $\psi(\cdot)$ is the solution of equation \eqref{eqGRA} with $\psi(0)=\phi_j(\gamma(0))$ corresponding to the control $\uu:[0,1/\eps]\rightarrow \omega$ defined by $\uu(t)=\gamma(\eps t)$, 
\[
p_1=|\cos \left( \bth(\al_+) -\bth(\al_-) \right) |, \quad
p_2=|\sin \left( \bth(\al_+) -\bth(\al_-) \right) |,
\]
and $\bth(\cdot)$ is defined as in Proposition \ref{tetha}. 
\ep

\proof 
We consider the Hamiltonian $H(\uu(t)), \ t \in [0,1/\eps]$.
Since the control function $\uu(\cdot)$ is not $\con^1$ at the singularity, we cannot directly apply the adiabatic theorem. Instead, we consider separately the evolution on the two subintervals (in time $t$) $[0,\tau_0 /\eps]$ and $[ \tau_0 /\eps,1/\eps]$. 

Since the eigenstates $\phi_j(\uu(t)),\phi_{j+1}(\uu(t))$ are piecewise $\con^1$ we can apply Theorem~\ref{effective} in order to study the evolution inside the space $P_{\uu(t)}\H$.
We can then construct the effective Hamiltonian, which is diagonal on both intervals (in time $\tau$) $[0,\tau_0]$ and $[\tau_0,1]$. 
Remark that the operator-valued function ${\cal U}(\cdot)$, defined in equation~\eqref{Urep}, has a discontinuity at $\tau_0 $ but has continuous extensions on both intervals $[0,\tau_0]$ and $[\tau_0,1]$. 

Let $\phi_j^{\pm}=\lim_{\tau \to \tau_0^{\pm}}\phi_j(\gamma(\tau))$. Integrating the effective Hamiltonian we get  
\[
U_a^{\eps}(\tau_0,0) \psi(0)= e^{i \varphi} \phi_j^-
\]
for some $\varphi \in \R$ . By Proposition \ref{tetha} we have
\begin{align*}
\phi_j^+&= \cos \left (\boldsymbol{\vartheta}(\al_+) - \boldsymbol{\vartheta} (\al_-) \right)  \phi_j^- + \sin \left (\boldsymbol{\vartheta}(\al_+) - \boldsymbol{\vartheta} (\al_-) \right)  \phi_{j+1}^-  \\
\phi_{j+1}^+&= -\sin \left (\boldsymbol{\vartheta}(\al_+) - \boldsymbol{\vartheta} (\al_-) \right)  \phi_j^- + \cos \left (\boldsymbol{\vartheta}(\al_+) - \boldsymbol{\vartheta} (\al_-) \right)  \phi_{j+1}^- . 
\end{align*}
Then, since the effective Hamiltonian is diagonal, we get
\[
U_a^{\eps}(1,0) \psi(0)= e^{i \vartheta_j} \cos \left (\boldsymbol{\vartheta}(\al_+) - \boldsymbol{\vartheta} (\al_-) \right) \phi_j(\gamma(0)) + e^{i \vartheta_{j+1}} \sin \left (\boldsymbol{\vartheta}(\al_+) - \boldsymbol{\vartheta} (\al_-) \right) \phi_{j+1}(\gamma(0)) ,
\]
and then, applying the adiabatic theorem, 
\[
\|\psi(T)- p_1e^{i\vartheta_j}\phi_j(\gamma(0))-p_2e^{i\vartheta_{j+1}}\phi_{j+1}(\gamma(0))\|\leq \widehat{C} \eps
\]
where $\widehat{C}$ is a constant depending on the gap and on $\gamma$. \hfill $\Box$


\brem \label{remark2level}
For control purposes, it is interesting to consider the case in which the initial probability is concentrated in the first level, the final occupation probabilities $p_1^2$ and $p_2^2$ are prescribed, and there is an integral curve of $\Xpao$ connecting $\uu^0$ to the singularity. Except for the special cases $p_1^2=0,1$, there are exactly two integral curves of $-\Xpao$ starting from the singularity that 
realize the required splitting (in the sense of Proposition~\ref{prop-2level}). 

Choosing $\beta \in [0,\pi/2]$ such that $(p_1,p_2)=(\cos \beta,\sin \beta)$, we obtain that the two possible values for $\al_+$ are
\[
\al_+=\bth^{-1}\left( \beta + \bth(\al_-) + k_+\pi \right) \qquad
\al_+=\bth^{-1}\left( -\beta + \bth(\al_-) + k_-\pi \right),
\]
where $k_+,k_-\in \mathbb{Z}$ are chosen in such a way that $\  \beta + \bth(\al_-) + k_+\pi\ $ and $ \ -\beta + \bth(\al_-) + k_-\pi \ $ belong to the range of $\bth$. 

If $(p_1^2,p_2^2)=(0,1)$, then the path is unique with $\al_+=\al_-+\pi$, while if $(p_1^2,p_2^2)=(1,0)$, then the unique path satisfies $\al_+=\al_-$.
\erem 

\noindent
{\bf Proof of Theorem~\ref{asc-paolo}.}
For simplicity, we consider the case in which $\psi(0)=\phi_0(\uu^0)$. The general case can be treated similarly.

Recall that for any conical intersection between two eigenvalues of a \mariospectrum\ there exists a neighborhood of the intersection where the two eigenvalues are well separated from the rest of the spectrum.
Let us consider these neighborhoods for the intersections $\bar \xx_j, j=0,\ldots,k-1 $, and let us call them $\omega_j$. Define on each $\omega_j \setminus \{\bar \xx_j \}$ the vector field $\Xpao^j$ as in Section \ref{una-a-caso}. 

We construct the path $\gamma(\cdot)$ as described below. 

First choose a smooth path $\sigma_0(\cdot)$ starting from $\uu^0$ and reaching $\omega_0$ along which all the eigenvalues in $\Sigma$ are simple. Concatenate $\sigma_0$ with an integral curve of 
$\Xpao^0$ that reaches the point $\bar \xx_0$.
Then choose $\al_+^0$ as one of the angles realizing, for the two-levels system associated with the energy levels $\lam_0,\lam_1$, the splitting from $(1,0)$ to $(p_1^2,1-p_1^2)$, as explained in Remark~\ref{remark2level}, and continue the path with the integral curve of $-\Xpao^0$ with outgoing tangent parallel to $(\cos \al_+^0,\sin\al_+^0)$.

Join the latter with a smooth path $\sigma_1(\cdot)$ connecting $\omega_0$ to $\omega_1$ along which all the eigenvalues in $\Sigma$ are simple, and then prolong it with an integral curve of  
$\Xpao^1$ that reaches the point $\bar \xx_1$. 
As above, compute an angle $\al_+^1$ that realizes the splitting (for the two-levels system associated with the energy levels $\lam_1,\lam_2$) from $(1-p_1^2,0)$ to $(p_2^2,1-(p_1^2+p_2^2))$, and, as above, continue the path with the integral curve of $-\Xpao^1$ with outgoing tangent parallel to $(\cos \al_+^1,\sin\al_+^1)$.

Repeat this procedure iteratively until the required spread is realized. Then reach the final point  $\uu^1$ with a path along which all the eigenvalues are simple. We assume without loss of generality that the final time is equal to one. 

For $\eps>0$, consider the Hamiltonian $H(\uu(t))=H(\gamma(\eps t))$, and set $\tau=\eps t$. 

As long as $\gamma(\tau) \in \mathbb{R}^2\setminus \cup_{i=0}^{k-1} \omega_i$, we approximate the dynamics of $H(\uu)$ using the adiabatic Hamiltonian 
\begin{equation} \label{adiatutto}
h_a(\tau) = H(\gamma(\tau)) -i\eps \sum_{l=0}^k  P_l(\tau) \dot{P}_l(\tau) -i\eps  P_{\Sigma}^{\bot}(\tau) \dot{P}_{\Sigma}^{\bot}(\tau)
\end{equation} 
where $P_l(\tau) $ is the spectral projector onto the eigenspace relative to $\lambda_l(\gamma(\tau))$
and $P_{\Sigma}^{\bot}(\tau)=\mathrm{id} - \sum_{l=0}^k  P_l(\tau) $. 

The evolution associated with \eqref{adiatutto} conserves the occupation probabilities relative to each energy level in $\Sigma$, therefore the evolution of $H(\gamma(\tau))$ approximately conserves these occupation probabilities, with an error of the order $\eps$, as prescribed by the adiabatic theorem (see Remark~\ref{nenciu}).

For $\gamma(\tau) \in  \omega_j, \ j=0,\ldots,k-1$, we use instead the adiabatic Hamiltonian
\beq  \label{adiaj}
h_a(\tau) = H(\gamma(\tau)) -i\eps P_{j,j+1}(\tau) \dot{P}_{j,j+1}(\tau) -i\eps \sum_{\substack{l=0 \\ l\neq j,j+1} }^k  P_l(\tau) \dot{P}_l(\tau) -i\eps  P_{\Sigma}^{\bot}(\tau) \dot{P}_{\Sigma}^{\bot}(\tau)
\eeq 
where $P_{j,j+1}(\tau)$ is the spectral projector relative to $\{\lambda_j(\gamma(\tau)),\lambda_{j+1}(\gamma(\tau))\}$. 

The evolution associated with \eqref{adiaj} conserves the occupation probabilities relative to the band
$\{\lambda_j,\lambda_{j+1}\}$, to any other energy level in $\Sigma$ and to its remainder in the spectrum. Moreover, thanks to the choice of the field $\Xpao^j$, we can also compute the evolution given by  \eqref{adiaj} inside the band
$\{\lambda_j,\lambda_{j+1}\}$ (which is the one described in Proposition~\ref{prop-2level}).

We end up with final state $ \psi(1/\eps)$ satisfying
\[\|\psi(1/\eps)-\sum_{l=0}^k p_l e^{i\vartheta_l} \phi_l(\uu^1)\|\leq C \eps,
\]
for some $\vartheta_0,\ldots,\vartheta_k\in \mathbb{R}$ and some $C>0$ determined by the adiabatic approximation. Thus the system is approximately spread controllable and the theorem is proved.
\hfill $\square$

\section{Mildly mixing curves}\label{mild}

In the previous section we constructed some special curves along which the effective Hamiltonian has a simple form, whose evolution is quite easy to predict.
In this section, we consider more general curves passing through the singularities.

We prove below
 a variation of  Proposition~\ref{prop-2level}, which generalizes to broken curves the result in 
 \cite[Corollary 2.5]{teufel}:  if we choose any piecewise  regular curve with a vertex at the conical singularity, then  
we obtain a distribution of probability between the two levels 
similar to the one described by
Proposition \ref{prop-2level}.
In this case, if the final time is $1/\eps$, the error is of order $\sqrt{\eps}$. 

Moreover, we prove that the integral curves of $\Xpao$ are not the only ones that realize the best accuracy (that is, an error which is of order $\eps$ for a final time equal to $1/\eps$): indeed, this can be obtained with any curve whose first and second derivatives at the singularity are the same as those of an integral curve of $\Xpao$. 

Let us consider a $\con^2$ curve $\gamma : [0,\tau_0] \to \omega$ such that $\gamma(\tau_0)=0$   corresponds to a conical intersection between $\lambda_j$ and $\lambda_{j+1}$, and $\dot{\gamma}(\tau_0^-)\neq 0$. Assume moreover that $\phi_j,\phi_{j+1}$ are $\con^2$ along $\gamma$ (recall that this is true for analytic curves). Let us consider the Hamiltonian $H(\gamma(\eps t)), t\in [0,\tau_0/\eps]$, and the adiabatic Hamiltonian \eqref{adiaj}.
Up to a factorization of the trace, the effective Hamiltonian reads
\[
\Heff(\tau) = \begin{pmatrix} a(\tau) & -i \eps b(\tau) \\ i \eps b(\tau) & -a(\tau) \end{pmatrix}
\]
where 
\[
a(\tau)=\frac{\lam_{j+1}(\gamma(\tau))-\lam_j(\gamma(\tau))}{2}, \qquad 
b(\tau) = \langle \phi_{j+1}(\gamma(\tau)),\dot{\phi}_j(\gamma(\tau)) \rangle,
\] 
and the dynamical system associated with $\Heff$ is
\beq \label{Hz}
i\eps\dot z = \Heff z
\eeq
where $z\in\C^2,\ |z|=1$ (as usual, here and below the dot indicates the derivative with respect to $\tau$). 

Note that the condition of conical intersection implies the existence of two positive constants $C_1,C_2$ such that, for $\tau$ close to $\tau_0$,  $-C_2 \leq \dot a(\tau) \leq -C_1$ (as a consequence $C_1 |\tau-\tau_0| \leq a(\tau) \leq C_2 |\tau-\tau_0|$).
As for $b(\tau)$, it is $\con^1$ by hypothesis.

We set  $D(\tau) = \left( \begin{smallmatrix} a(\tau)/\eps & 0 \\ 0 & -a(\tau)/\eps\end{smallmatrix}\right)$, $U(\tau) = \exp \left(-i \int_0^{\tau} D(s) ds \right)$,
and we perform the change of variable $\zeta=Uz$, so that $\zeta$ evolves according to the dynamical system $i \dot{\zeta} = \widehat{\Heff}(\tau) \zeta$, where
\[
\widehat{\Heff}(\tau) = \frac{1}{\varepsilon}U(\tau)\Heff(\tau) U(\tau)^{-1} + \dot{U}(\tau)U(\tau)^{-1}=\ \left(\ba{cc} 0 & -i b(\tau)e^{\frac{2i}{\eps}\int_0^{\tau} a(s)ds} \\ i b(\tau)e^{-\frac{2i}{\eps}\int_0^{\tau} a(s)ds} & 0\ea\right)\,.
\]
Let us express the evolution operator for $\widehat{\Heff}$ in the form 
\[
M_{\eps}(\tau,0) = \left(\ba{cc} \nu(\tau) & \mu^*(\tau)e^{\frac{2i}{\eps}\int_0^{\tau} a(s)ds} \\ -\mu(\tau)e^{-\frac{2i}{\eps}\int_0^{\tau} a(s)ds} & \nu^*(\tau)\ea\right)\,.
\]
We claim that $\|M_{\eps}(\tau,0)-\mathrm{id}\| \leq C \sqrt{\eps}$, for some $C>0$.

From $i\frac{d}{d\tau}M_{\eps}(\tau,0)=\widehat{\Heff}(\tau) M_{\eps}(\tau,0)$ we get the equations
\[
\left\{\ba{l}
\dot \nu = \mu b\\
\dot \mu = \frac{2i}{\eps}\mu a - \nu b
\ea\right.
\]
with initial data $\nu(0)=1,\mu(0)=0$.
Since $\det M_{\eps}=1$ we have that $|\nu|$ and $|\mu|$ are bounded, and then, from the boundedness of $b$, we get that also $|\dot\nu|$ is bounded.

We recall that along the integral curves of $\Xpao$ the effective Hamiltonian is diagonal and its evolution is exactly $U(\tau)$, so that the equations above are solved by $\nu\equiv 1$, $\mu \equiv 0$. If $b(\cdot)$ is not identically equal to zero, then the evolution is not exactly diagonal, but it mixes the two components of $\zeta$ (and $z$). The error done by assuming the evolution diagonal can be estimated by evaluating the term $\mu$.

By variation of constants we have
\[\mu(\tau) = - \int_0^{\tau} e^{\frac{2i}{\eps} \int_s^{\tau} a(r)dr} b(s)\nu(s) ds.\] 
We can rewrite 
\begin{align}
\mu(\tau_0) &= \int_0^{\tau_0} -a(s) \frac{e^{\frac{2i}{\eps} \int_s^{\tau_0} a(r)dr} b(s)\nu(s)}{a(s)} \;ds\nn\\
&= \Big[\frac{\eps}{2i} \left( e^{\frac{2i}{\eps} \int_s^{\tau_0} a(r)dr}-1 \right) \frac{b(s)\nu(s)}{a(s)}\Big]_0^{\tau_0} \!\!-\!\! \int_0^{\tau_0} \frac{\eps}{2i}\left(e^{\frac{2i}{\eps} \int_s^{\tau_0} a(r)dr}-1 \right) \frac{d}{ds} \left( \frac{b(s)\nu(s)}{a(s)} \right) \; ds\nn\\
&= \frac{\eps}{2i} b(\tau_0)\nu(\tau_0)\lim_{\tau\to \tau_0}  \frac{e^{\frac{2i}{\eps} \int_{\tau}^{\tau_0} a(r)dr}-1}{a(\tau)}  - \frac{\eps}{2i} \left( e^{\frac{2i}{\eps} \int_0^{\tau_0} a(r)dr}-1 \right) \frac{b(0)}{a(0)} + \nonumber \\ 
&-\frac{\eps}{2i} \int_0^{\tau_0} \left(e^{\frac{2i}{\eps} \int_s^{\tau_0} a(r)dr}-1 \right) \frac{d}{ds} \left( \frac{b(s)\nu(s)}{a(s)} \right) \; ds \label{perparti}.
\end{align}
Since $\dot{a}(\tau_0)\neq 0$, we obtain that $\lim_{t\to \tau_0}  \frac{e^{\frac{2i}{\eps} \int_t^{\tau_0} a(r)dr}-1}{a(\tau)}=0$.
The second term in the equation above is of order $\eps$. Then, we are left to estimate the integral term. We have
\bqn
\left|\int_0^{\tau_0} \left(e^{\frac{2i}{\eps} \int_s^{\tau_0} a(r)dr}-1 \right)\frac{d}{ds} \left(\frac{b(s)\nu(s)}{a(s)}\right) \; ds\right| \leq C \int_0^{\tau_0} \frac{|e^{\frac{2i}{\eps} \int_s^{\tau_0} a(r)dr}-1|}{a^2(s)} \; ds. \nn\label{stima-easy}
\eqn
We consider the change of variables $s\mapsto \xi_{\eps}(s)=\frac{2}{\eps} \int_s^{\tau_0} a(r)dr$ so that 
\[
\int_0^{\tau_0} \frac{|e^{\frac{2i}{\eps} \int_s^{\tau_0} a(r)dr}-1|}{a^2(s)} ds = \frac{\eps}2 \int_0^{\xi_{\eps}(0)} \frac{|e^{ix}-1|}{a^3(\xi_{\eps}^{-1}(x))} \;dx.
\]
From the estimates above on $a$ we easily get that $a^3(\xi_{\eps}^{-1}(x))\geq \bar C \eps^{3/2} x^{3/2}$ for a suitable positive constant $\bar C$.

Since $\int_0^{+\infty}  \frac{|e^{iz}-1|}{z^{3/2}} dz\leq \int_0^{+\infty}\min \{z^{-1/2},2z^{-3/2}\}dz< +\infty$, we immediately obtain that  the integral in \eqref{perparti} is of order $\sqrt{\eps}$.
Therefore, $\|M_{\eps}(\tau_0,0) - \textrm{id}\|$ is of order $\sqrt{\eps}$.

If $\gamma$ is defined also for $\tau>\tau_0$ and is globally $\con^2$,  we recover Corollary 2.5 in \cite{teufel}.
If, instead, $\gamma$ is 
continuous and piecewise $\con^2$, with different tangent directions at the singularity, then we can repeat the same argument as in Proposition~\ref{prop-2level}: at the singularity the limit basis rotates instantaneously and we consider separately  the evolution of two different adiabatic Hamiltonians. The rotation of the limit basis spreads the probabilities as described by equations  \eqref{transangolo1}-\eqref{transangolo2},  and this leads to a controllability result in the spirit of Theorem~\ref{asc-paolo}, where the error is of order $\sqrt{\eps}$ if the final time is $1/\eps$. 

The following result shows that the value of $b(\tau)$ at the instant where the curve attains the singularity depends only on the 2-jet of the curve at the singularity. This allows us, using piecewise analytic curves that have the same 2-jet at the singularity as an integral curve of $\Xpao$, to obtain a controllability result equivalent to Theorem~\ref{asc-paolo} (see Proposition~\ref{last-p}). 

\bl \label{2-jet}
Let $\gamma(\cdot)$ and $\tg(\cdot)$ be two $\con^2$ curves on $\omega$ such that $\gamma(\tau_0)=\tg(\tau_0)=0$, where $0$ is a conical intersection between $\lambda_j$ and $\lambda_{j+1}$, with $\dot{\gamma}(\tau_0)=\dot{\tg}(\tau_0)\ne 0$ and $\ddot{\gamma}(\tau_0)=\ddot{\tg}(\tau_0)$.
Let $\eta(\tau) = \langle \phi_{j+1}(\gamma(\tau)),\dot{\phi}_j(\gamma(\tau))\rangle$ and $\widehat{\eta}(\tau)= \langle \phi_{j+1}(\tg(\tau)),\dot{\phi}_j(\tg(\tau)) \rangle$.
Then $\lim_{\tau\rightarrow \tau_0} ( \eta(\tau) -\widehat{\eta}(\tau) ) =0$.
\el

\proof
First of all, we remark that $|\gamma(\tau)-\tg(\tau)|=o((\tau-\tau_0)^2)$ and $|\dot{\gamma}(\tau)-\dot{\tg}(\tau)|=o(|\tau-\tau_0|)$. 

As in the proof of Proposition \ref{limitati}, we can prove that there exists $C>0$ such that for every point $\uu\in \omega \setminus \{0\}$ and any unitary vector $\mathbf{w}\in \mathbb{R}^2$  the directional derivative along $\bf{w}$ satisfies $\|\partial_{\bf w} \phi_l(\uu)\|\leq C/|\uu|$. Then we obtain that $\|\phi_l(\gamma(\tau))-\phi_l(\tg(\tau))\| = o(|\tau-\tau_0|), \ l=j,j+1$.
Moreover, by \eqref{obs1} we know that the eigenvalues are Lipschitz in a neighborhood of the intersection.

From \eqref{stima} we have
\begin{align*}
\eta(\tau) - \widehat{\eta}(\tau) &= \frac{\langle \left( \phi_j(\gamma(\tau)) -\phi_j(\tg(\tau))\right) , \left(\dot{\gamma}_1H_1+\dot{\gamma}_2 H_2\right) \phi_{j+1}(\gamma(\tau)) \rangle}{\lambda_{j+1}(\gamma(\tau))-\lambda_{j}(\gamma(\tau))}  \\
&+\frac{\langle \phi_j(\tg(\tau)) , ((\dot{\gamma}_1-\dot{\tg}_1)H_1+(\dot{\gamma}_2-\dot{\tg}_2) H_2) \phi_{j+1}(\gamma(\tau)) \rangle}{\lambda_{j+1}(\gamma(\tau))-\lambda_{j}(\gamma(\tau))}\\
&+\frac{\langle  \phi_j(\tg(\tau)) , ( \dot{\tg}_1 H_1+\dot{\tg}_2 H_2 ) \left( \phi_{j+1}(\gamma(\tau)) -\phi_{j+1}(\tg(\tau))\right) \rangle }{\lambda_{j+1}(\gamma(\tau))-\lambda_{j}(\gamma(\tau))}\\
&+\langle  \phi_j(\tg(\tau)) , (\dot{\tg}_1H_1+\dot{\tg}_2 H_2) \phi_{j+1}(\tg(\tau)) \rangle \left(\frac{1}{\lambda_{j+1}(\gamma(\tau))-\lambda_{j}(\gamma(\tau))} -\frac{1}{\lambda_{j+1}(\tg(\tau))-\lambda_{j}(\tg(\tau))} \right).
\end{align*}

By previous estimates it 
follows that all the terms in the right-hand side of the equation above go to zero as $\tau$ goes to $\tau_0$.
\eproof

\bp\label{last-p}
Let $\uu=0$ be a conical intersection between the eigenvalues $\lam_j,\lam_{j+1}$ and let $\phi_j^0,\phi_{j+1}^0$ be limits as $\tau \to 0^+$ of the eigenstates $\phi_j(r(\tau) ),\phi_{j+1}(r(\tau) )$, respectively, for $r(\tau)=(\tau,0)$. Let $\gamma:[0,1]\rightarrow \omega$ be a curve such that there exists $\tau_0\in (0,1)$ with $\gamma(\tau)=0$ if and only if $\tau=\tau_0$, $\gamma$ analytic on $[0,\tau_0]$ and $[\tau_0,1]$, and $\dot{\gamma}(\tau_0^{\pm})\neq 0$. Let $\al_-$ and $\al_+$ be the angles describing respectively the inward and the outward tangent direction at the singularity, as in \eqref{eq-2level}.
Assume that the integral curves of $\Xpao$ having the same inward and outward tangents as $\gamma$ at the singularity possess also the same 2-jet as $\gamma$ at the singularity.
Then there exists $C>0$ such that for any  $\varepsilon >0$
\bqn 
\|\psi(1/\eps)- p_1e^{i\vartheta_j}\phi_j(\gamma(0))-p_2e^{i\vartheta_{j+1}}\phi_{j+1}(\gamma(0))\|\leq C \eps
\eqn 
where $\vartheta_j,\ \vartheta_{j+1} \in \R$,  $\psi(\cdot)$ is the solution of equation \eqref{eqGRA} with $\psi(0)=\phi_j(\gamma(0))$ corresponding to the control $\uu:[0,1/\eps]\rightarrow \omega$ defined by $\uu(t)=\gamma(\eps t)$, and 
\[
p_1=|\cos \left( \boldsymbol{\vartheta}(\al_+) -\boldsymbol{\vartheta}(\al_-) \right) | \quad
p_2=|\sin \left( \boldsymbol{\vartheta}(\al_+) -\boldsymbol{\vartheta}(\al_-) \right) |,
\]
with $\boldsymbol{\vartheta}(\cdot)$ defined as in Proposition \ref{tetha}. 
\ep
\proof
By Lemma~\ref{2-jet} 
the function $b(\tau)=\langle \phi_{j+1}(\gamma(\tau)),\dot{\phi}_j(\gamma(\tau)) \rangle$ goes to zero as $\tau$ goes to $\tau_0$. Moreover, the analyticity of $\gamma(\cdot)$ easily
implies that the term $\frac{d}{ds} \big(b(s)\nu(s)/a(s)\big)$ appearing in \eqref{perparti} is bounded. Thus $|\mu(\tau_0)|\leq C\eps$ for a suitable $C>0$.

Then 
$\|M_{\eps}(\tau_0,0) - \textrm{id}\|$ is of order ${\eps}$, where $M_{\eps}$ is the evolution operator defined above. We can obtain an analogous estimate for $\|M_{\eps}(\tau,\tau_0) - \textrm{id}\|,\ \tau >\tau_0$. This completes the proof.
\eproof

\bibliographystyle{plain}
\bibliography{biblio-conica-TAC}

\end{document}